\documentclass[12pt,a4paper]{article}
\input{epsf}

\usepackage{amsmath,amssymb,amsfonts,graphicx}
\usepackage{amsmath,amssymb,amsfonts,graphicx,amsthm}
\usepackage{cite}

\usepackage{tikz-cd} 
\usepackage{stix}

\usepackage[T1]{fontenc}
\usepackage{amsmath,amssymb,amsfonts,graphicx,amsthm,mathrsfs}
\usepackage{cite}
\usepackage[T1]{fontenc}
\usepackage[utf8]{inputenc}
\usepackage{authblk}
\usepackage[colorlinks]{hyperref}
 
\newtheorem{theo}{Theorem}
\newtheorem{lemma}{Lemma}

\newtheorem{cor}{Corollary}

\newtheorem{rem}{Remark}
\newtheorem{defi}{Definition}
 \newcommand{\R}{\mathbb{R}}
 \newcommand{\C}{\mathbb{C}}

\usepackage{enumerate}
\allowdisplaybreaks
\usepackage{hyperref}
\usepackage{stackrel}

\usepackage{array}

\date{}
\begin{document}
 
\title{On a Class of  Singular Complex Manifolds}
\author[1]{Alireza Bahraini\thanks{bahraini@sharif.edu}}
%\author[1]{Saeed Sadeghi\thanks{saeedsadeghi91@gmail.com}}
%\author[2]{Ali Mohammad-Djafari\thanks{djafari@ieee.org}}
\affil[1]{ Department of Mathematical Sciences, Sharif University of Technology, P.O.Box 11155-9415, Tehran, Iran.}
%\affil[2]{ Laboratoire des Signaux et Systèmes (L2S),  (CNRS-SUPELEC-UPS 11) Sup\'elec, Plateau de Moulon, 91192 Gif-sur-Yvette France}
\renewcommand*{\Authands}{, }
\maketitle

\begin{abstract} We introduce a new class of singular complex manifolds and   we develop  a
degenerate Kodaira-Hodge theory   for this class of singular  manifolds. 
\end{abstract}
%%%
%\maketitle
%%
%\centerline{\Large{$\S$ Introduction}}\\

\vspace{0.1cm}


\vspace{1.5cm}
\section{Introduction}

Recently we have  proved that the solutions to   degenerate complex  Monge-Amp\`ere equations (DCMA)  have a regular behavior  along the locus of  degeneration   which is supposed to be a smooth divisor  in a K\"ahler manifold \cite{bah}.    More precisely given a K\"ahler manifold $(X,\omega )$ of dimension $n$ and a smooth divisor $[D]$ in $X$ then the solution to the  following  degenerate  complex Monge Amp\`ere equation
\[
(\omega+\partial\bar{\partial} \phi)^n=e^F |S|^2 \omega^n
\]

is  $C^{\infty}$ in a neighborhood of  $D$ if the function $F$ belongs to $C^{\infty} (X)$. Here $S\in H^0 (X, [D])$ is a holomorphic section of $[D]$ vanishing along $D$.

If $D$ is assumed to be a canonical divisor in a complex surface of general type denoted by $X$, then one can apply  the solutions to  the above Monge-Amp\`ere  equation to  equip   $X\setminus D$, with  a hyperk\"ahler structure whose  behavior  along $D$,  by the above mentioned  thereom on DCMA equations is now clearly understood.
If $\mathbb{S}$ denotes the twistor sphere  on  $X\setminus D$, then it can be shown that  $\mathbb{S}$  admits a smooth extension to  the whole space $X$.

All the complex structures in $\mathbb{S}$  except the initial one are singular along $D$. The type of singularities  are known as   involutive structures following Treves \cite{tre}, where   $D$ constitute  the corresponding  characteristic set. The symplectic counterpart  has also been studied 
by Guillemine et al. \cite{GM}.
   Our goal in this paper is  to study   the singular complex manifold obtained by a degenerate complex structure $J'$ in $\mathbb{S}$,  which is orthogonal to the initial complex structure $J$ on $X$.  
It turns out  that, some of the most fundamental theories about K\"ahler manifolds can    be proven for this new class of singular complex manifolds.

Firstly $\bar{\partial}$ solvability in semi-local Stein type  domains  containing the  type of singularities like $J'$ along a divisor inside the domain is established. Treves et al. had already studies the  $\bar{\partial}$ solvability problem for the hypoanalytic manifolds with one dimensional  characteristic sets. In those cases the solvability does not hold true.  We apply the solvability along with a  weak  $L^2$-Hodge theory  to derive a string Hodge theory as well as a Kodaira  theory for $(X, J')$. In fact $L^2$ Hodge theory  for manifolds with non-isolated conical singularitie  has already been discussed in several papers by Cheeger and other co-authors. \cite{cheeger3} \cite{cheeger4}.  The singular complex manifolds  gthat arise from the solutions to DCMA have the same type of singularities as thosse studied there.  The difference is that we have singularity over complex structures too.  For this reason we have introduced in 
definition 2 a generalized class of singular complex manifolds for whcih the weak  $L^2$ Hodge theory can hold true.

   The initial motivation for this work  was to generate lagrangian cycles which represent the Poincar\'e  dual of certain cohomology classes, orthogonal to the first Chern class of a given  K\"ahler surface of general type. 
 According to mirror symmetry conjecture \cite{ktc} lagrangian submnaifolds
are of the same importance as  the holomorphic objects for dual Calabi-Yau manifolds. Although the realization of  lagrangian submanifolds is not straightforward due to their  nonlinear equations.  

In order to develop Kodaira-Hodge theory we start by developing an $L^2$ Hodge theory for singular complex manifolds introduced in definition (\ref{def2}). The  real structure of the singularities we consider is closely related to  the structure  studied by Cheeger  and Dai \cite{cheeger3} while the singularities they consider is sole;y over the metric  and not over the complex structure.

We start by describing in section \ref{sec2}, the degenerate hyperK\"ahler structure obtained through the solutions to DCMA equations on a  K\"ahler surface containing a smooth canonical divisor.  In section \ref{sec3} we strudy the dergenerate complex structure $J'$ orthogonal to the initial complex structure $J$ in the corresponding twistor sphere and we provide more examples  that possess  the same type of  singularity.

%In section \ref{sec41} we recall the $L^2$ Hodge theory developed by Cheeger and Dai \cite{cheeger3} which is closely related with   the singularities of the metric  studied in the current  paper.  In the definition \ref{def2} of  subsection \ref{sec42} we introduce singular complex manifolds for which we can prove an $L^2$ Hodge theory. 

In sections \ref{weakh} and \ref{sec6}  weak Hodge theory for the  class of singular manifolds obtained through degenerate hyperk\"ahler structure  of section \ref{sec2}  and the isomorphism with the  complex obtained by the closure $\bar{\partial }^c _{J'}$  %and  the isomorphism with $L^2$-harmonic forms 
are established. As a result   a  strong  Hodge theory for our  singular complex manifolds is deduced.. 

In section \ref{sec8}, we develop  a Kodaira embedding theorem for the complement of an arbitrarily small  neighborhood $U_{\epsilon}$ of $D$ in $X$.
As will be discussed, the degenerate complex structure $J'$  induces a pseudo-concave structure on  $X\setminus U_{\epsilon}$.
It is well known that in complex dimension $2$ a generic  CR manifold (here $\partial U_{\epsilon}$) does not admit any embedding to a projective space.  Our result shows in particular  that the degenerate complex structures in the hyperk\"ahler structure  are   exceptional as can be expected.  
     %These are in fact   very general singularities  that have been first introduced and studied by  F.Treves, P.D. Cordaro, M.G. Eastwood  %and several other mathematicians generally in the
 %framework of overdetermined systems of PDE
  %(See \cite{tre}).
\\

 We establish in section \ref{solv} a $\bar{\partial}_{J'}$-solvability for the structure $J'$ in  a neighborhood of $D$. This  basic result  is  applied to prove the fundamental result on  weak  for Hodge   theory  studied in sections \ref{weakh}.
 \\

\section{DCMA equations and Degenerate HyperK\"ahler structure}\label{sec2}
  \label{dgma}

Consider a K\"ahler manifold $(X, \omega, D)$ where $\omega$ is a K\"ahler metric  and $D\subset X$  denotes a  smooth divisor on $X$.
  In (\cite{bah})  we have proved that there exists a $(1,1)$-form $\omega'\in \Omega ^{1,1} (X)$ with the following properties

i)  $\omega'$ has the same cohomology class as $\omega'$

\[
[\omega']=[\omega]\in H^{1,1} (X).
\]

ii) $\omega'$ is a K\"ahler form on $X\setminus D$.

iii) $j_D^* \omega'$ is a K\"ahler metric on $D$ and $\omega'$ degenerates transversally along $D$.  More precisely there exists a normal bundle $N_D\subset T' _{\mathbb{C}} X|_{D}$  such that $\omega'$ degenerates along the direction $N_D$, and
\[
\omega' \wedge \omega' \sim O(|S|^2)
\] 
  near $D$, 
where $S$ is the holomorphic section of $[D]$  vanishing on $D$.

We have then shown that the following DCMA equation 
\begin{equation}\label{mgap}
(\omega+\partial\bar{\partial}\phi)^2 = e^G  |S|^2 \omega'^2,
\end{equation}

admits a smooth solution $\phi\in C^{\infty} (X)$ if $G\in C^{\infty} (X)$ satisfies the required volume condition
$\int_X e^F \omega'^2 = \int_X  \omega'^2$.

Now  if   $D$  is a smooth canonical divisor in  the complex surface $X$ and if $\Omega$ denotes a holomorphic $(2,0)$-form
  vanishing along $D$, then   by solving the degenerate Monge-Amp\`ere equation
\begin{equation}\label{degm}
\omega'\wedge \omega' = \lambda \Omega\wedge \bar{\Omega} 
\end{equation}

for appropriate constant $\lambda$, we obtain a Ricci-flat metric $g'$ on $X\setminus D$ degenerating  transvarsally along $D$.

%\textbf{Degenerate quaternionic structre on $X$ }. 
Let 
\[
\Omega = \alpha+i\beta
\]

be the decomposition of $\Omega$ into real and imaginary parts.  Consider   the map

\begin{align*}
&p:S^2\rightarrow \Omega ^2 (X),\notag \\ 
 &(a,b,c)\rightarrow a\omega' +b \alpha +c \beta,
\end{align*} 

 where $S^2=\{(x,y,z)| x^2+y^2+z^2=1\}\subset \mathbb{R}^3$ is the unit sphere.
Then, for any positively oriented  orthonormal frame $(v_1,v_2, v_3)$ of unit vectors $v_1,v_2, v_3\in S^2$  the  subbundle

\begin{equation}\label{invo1}
\mathcal{V}:=Ker \left(p\left(v_1\right)+ip\left(v_2\right)\right) \subset T_{\mathbb{C}}  X,
\end{equation}
 
defines an involutive structure on $X$ (\cite{tre}).  This means that the complex subbundle  of $T_{\mathbb{C}}X$ defined by  $\mathcal{V} $  satisfies Frobenious integrabiluty condition $[\mathcal{V}, \mathcal{V}]\subset  \mathcal{V}$. Moreover we have $\mathcal{V}\cap \bar{\mathcal{V}} =\{0\}$ and $\mathcal{V}\oplus \bar{\mathcal{V}} =T_{\mathbb{C}}$ on $X\setminus D$. In other words   we have a sphere $\mathbb{S}$ of involutive structures on $X$  inducing  a hyperK\"ahler structure on $X\setminus D$.
%In this paper we  study  those degenerate complex structures in the twistor sphere $\mathbb{S}$, which are perpendicular to the initial %complex structure $J$ on $X$.   We start by describing explicitly in simple examples the local behavior of these degenerate complex %structures.
  \\

\textbf{Notation} The degenerate  complex structure compatible with  the degenerate K\"ahler form $\Re{\Omega}$   and the degenerate Ricci Flat metric $g'$ is denoted by $J'$. The manifold $(X,g', J')$ constitute the main class of examples we study in this paper.
   
\section{Involutive Structures and Hypoanalytic manifolds following F.Tr\`eves\cite{tre}}\label{sec3}
  
 We recall some definitions from (\cite{tre}) regarding Hypoanalytic manifolds and involutive structures.
 Let $M$ be a
 smooth manifold. The notion of \emph{involutive} structure on $M$ consists of
  a vector subbundle $\mathcal{V}$ of the complexified tangent
  bundle $T_{\mathbb{C}} M $ satisfying the Frobenious condition
  $[\mathcal{V}, \mathcal{V} ] \subset \mathcal{V}$. The structure is locally integrable if
  the subbundle $T'$ of $T_{\mathbb{C}} ^* M $ orthogonal to $\mathcal{V}$
   for the natural duality between tangent and cotangent vectors
   is locally spanned by exact differentials.
   The
  \emph{charateristic} set of the structure is the subset $T^ {\circ}$  of
the real tangent bundle $T^* M $  satisfying $T^{\circ} \otimes_{\mathbb{R}} \mathbb{C}=T'\cap \bar{T}'$. It is also defined  as to be  the common zeros
of the symbols of the sections of $\mathcal{V}$ i.e. $T^{\circ}
= T^* M \cap T'$.

A smooth manifold $M$, if it is equipped with a locally
integrable structure, can be covered by open sets  each of which
is mapped into $\mathbb{C}^m $ by a $C^\infty$ map
$Z=(Z_1,...,Z_m ) $. A hypo-analytic structure on $M$ is defined
as follows:

\begin{defi}\label{def0}
Let $M$ be a $C^\infty$ manifold. By a hypo-analytic structure on
$M$  we mean  the data of an open covering $\{\Omega_i\}$ of $M$
and for each index $i$, of a $C^\infty$ map
$Z_i=(Z_{i,1},...,Z_{i,m}): \Omega_i \rightarrow \C^m$ with $m\geq
1$ independent of $i$, such that the following is true:

1) $dZ_{i,1},...,dZ_{i,m}$ are $\C$-linearly independent at
each point;

2) if $i\neq k$ for each point $p$ of $\Omega_i \cap \Omega_k$
there is a holomorphic map $F^i_{k,p}$ of an open neighborhood of
$Z_i (p)$ in $\C^m$ into $\C^m$ such that $Z_k=F_{k,p}^i o Z_i$
in a neighborhood of $p$ in $\Omega_i \cap \Omega_k$.
\end{defi}
 Most of the examples studied by Tr\`eves et al. have characteristic
 sets of dimension equal to one.
The following example of involutive structures could be
considered as a generalization
 of Mizohata structures (see e.g. \cite{tre} ) to complex dimension greater
 than one.
\\

The following example of involutive structures could be
considered as a generalization
 of Mizohata structures (see e.g. \cite{tre} ) to the case where characteristic
 set is given  by a  submanifold of dimension equal to two.  
\\
\\

 \emph{Simple branched covering.}
 Let
\[
\Omega \subset \mathbb{C}^2=\{(w,z)|w=u+iv,z=x+iy)\}
\]
 be an open subset and let $
 \mathbb{H}=\{(U+iV,X+iY)| U,V,X,Y\in \R\}$ be equipped with its
 standard quaternionic structure $\{I_0,J_0\}$ and let $\pi: \C ^2 \rightarrow \mathbb{H}$ be
 the application defined by
 \[
\pi(w,z)=(w,z^2).
 \]
 Assume that the complex coordinates associated to $J$ is given by $(U+iX,V-iY)$. The restriction of $\pi$ to $\Omega$,
\[
\pi|_{\Omega}:\Omega \rightarrow  \mathcal{W},
\]
  creates  a double covering of  an open subset $\mathcal{W}$ of
  $\mathbb{H}$, ramified  along
$\R^2\times\{0\} \cap \Omega $. Thus,
 $\pi$ defines a hypo-analytic structure on
$\Omega$,  which is given by the following
holomorphic coordinates:
\[
\begin{matrix}
\tau=u+i(x^2-y^2)& \rho=v-i(2xy).
\end{matrix}
\]
The pullback $J' _0=\pi^*(J_0)$ of the complex structure $J_0$ defines 
a degenerate complex structure on $\Omega$, which is exactly of the same type we study in this paper. If $r
  ^2=x^2+y^2$ then $J'_0$ can be described as follows:

\begin{align*}
J'_0(dx)&=\frac{1}{2 r^2} (xdu+ydv)\\
J'_0(dy)&=\frac{1}{2 r^2}(xdv-ydu)\\
J'_0(du)&=-2xdx+2ydy\\
J'_0(dv)&=2xdy+2ydx.
\end{align*}

Following Andr\'e Weil  % \cite{weil}
 we have:
\begin{equation}\label{dbar}
\bar\partial_{J'_0} = \frac{1}{2}(d+i(J'_0)^{-1}dJ_0)
\end{equation}
which yields the following degenerate Cauchy-Reimann equations:

\[
2\rho ^2
\begin{pmatrix}
\partial_u \phi\\
\partial_v \phi
\end{pmatrix}
+i
\begin{pmatrix}
x & -y \\
y &  x
\end{pmatrix}
\begin{pmatrix}
\partial_x \phi\\
\partial_y \phi
\end{pmatrix}
=0
\]
or

\[
\begin{pmatrix}
\partial_x \phi \\
\partial_y \phi
\end{pmatrix}
= 2i
\begin{pmatrix}
x & y \\
-y & x
\end{pmatrix}
\begin{pmatrix}
\partial_u \phi \\
\partial_v \phi
\end{pmatrix}
\]

 Let  $A^k (\Omega)$ denote  the space of $C^\infty$ sections of the
vector bundle $\pi^* (\Lambda ^k(T^* U))$,  and  we define

\[
A^{p,q}(\Omega):=A^{p+q}(\Omega)\cap A^{p,q}(\Omega\setminus
\R^2\times \{0\}).
\]

In local coordinates a basis for $A^{1,0}$ is given by

\begin{equation}\label{wprim}
\begin{matrix}
\omega'_1=du+i(2xdx-2ydy), & \omega'_2=dv-i(2xdy+2ydx).
\end{matrix}
\end{equation}

  The dual
 basis $\{v'_1, v'_2 \}$  of the basis $\{\omega' _1, \omega ' _2 \}$, with respect to the standard metric,
    is equal to
\begin{equation}\label{v}
\begin{matrix}
v'_1=\frac{1}{2}(\partial_u+\frac{1}{2i}\frac{x\partial_x-y\partial_y}{x^2+y^2}),&
v'_2=\frac{1}{2}(\partial_v
-\frac{1}{2i}\frac{y\partial_x+x\partial_y}{x^2+y^2}).
\end{matrix}
\end{equation}

Hence we can define the operator $\bar{\partial}_{J'}$ as follows
\[
\bar{\partial}_{J'}:= v'_1 \omega' _1+ v'_2 \omega'_2.
\]

\begin{rem}\label{asasi}
The operator $L=\rho ^2
\bar\partial _{J'}$ discussed   by F.Treves in \cite{tre}, can be identified with the operator
$\bar\partial_{J'}$ defined by (\ref{dbar}) if the latter is considered in the weak sense with respect to the volume form of the degenerate metric.
\[
T_{\bar\partial _{J'} \alpha} (\phi)= \int \alpha \wedge
\bar\partial _{J'} \phi
\]
\end{rem}

There exists a straightforward generalization of  the above example  to
higher dimensional  complex spaces. Let
\[
\begin{matrix}
\pi^k:=\underbrace{\pi\times \pi \times ...\times \pi} \\

\hspace{1cm}\text{ k times}
\end{matrix}
\]
 and  define the application $\pi_{n,k}$ by:
\begin{equation}\label{pi}
\pi_{n,k}=\pi^{2k}\times Id_{n-2k}:\C^{2k}\times
\C^{n-2k}\rightarrow \mathbb{H}^k \times \C^{n-2k}
\end{equation}
Assume that $\mathbb {H}^k$ is equipped with the complex structure
$J^{2k}$  and consider the ordinary complex structure on
$\C^{n-2k}$ denoted by $I^{n-2k}$. It can be easily seen  that
$\C^{2n}$ with the chart defined by $\pi_{n,k}$ is a
hypo-analytic manifold and the corresponding complex structure
$J'_{n,k}=\pi_{n,k}^* (J^{2k}\times I^{n-2k})$ degenerates
transversally along $\R^{2k}\times 0$.
\\ 
\\

 \textbf{Global construction  of hypoanalytic manifolds.} Let $Y$ be a complex manifold of complex dimension
$2k$ with trivial canonical bundle equipped with a Ricci flat K\"ahler metric. Let $\Gamma\subset Y$ be a special
Lagrangian sub-manifold in $Y$. Assume that $u: (X,D)\rightarrow (Y,\Gamma)$ is a double covering ramified
along $D$ where $X$ is an analytic manifold and $D\subset X$ is a sub-manifold. Assume that for every
point $p\in D$ there exits local charts $\eta_1: U_p \rightarrow \C^{2k} $ and $\eta_2 :V_{u(p)}\rightarrow
\mathbb{H}^k$ such that $\eta_2\circ u \circ \eta_1 ^{-1}$ coincides with $\pi_{2k,k}$, defined by (\ref{pi}),
upto order 2 along $\eta_1(U_p \cap D)$. The map $u$ induces a hypo-analytic structure on $X$ whose
characteristic set is equal to $D$.
\\

\textbf{Main Example.} Our main class of hypoanalytic manifolds we are going to study consists of the  manifolds equipped with the involutiive structures $\mathcal{V}$  obtained through the 
solutions to the degenerate complex Monge Amp\`ere equations as described by relations (\ref{invo1}) in the subsection (\ref{dgma}) 
\\
\\

%***************************************************************************************************

\section { Hodge Theory on Degenerate Complex Manifolds}\label{sec4}

\subsection{$L^2$  Hodge Theory on manifolds with non-isolated conical singularities following Cheeger-Dai }\label{sec41}
 
The study of $L^2$ cohomology  of spaces with non-isolated conical singularities has been carried out  by Cheeger and Dai in \cite{cheeger3}.The local model of singularities they consider is described as follows. Given a fibration
\[
Z^n\rightarrow M^m\xrightarrow{\pi} B^l,
\]

of closed oriented smooth manifolds, if $C_{\pi}M$ denotes the mapping cylinder of $\pi$, then one can attach a cone to each of the fibers to get
\[
C_{[0,1]}(Z)\rightarrow C_{\pi}M\rightarrow B.
\]
Here $C_{[0,1]}(Z)=(0,1)\times Z$.
The submersion  metric 

\[
g_1= dr^2 +\pi^* g_B + r^2 g_Z
\]

induces a metric on the nonsingular part of $C_{\pi}M$, which extends to a singular metric on the whole space. The manifolds studied by Cheeger and Dai have the following form
\[
X=X_0\cup X_1\cup...\cup X_k,
\]

where $X_0$ is a compact smooth manifold with boundary, and  each $X_i$  ($i\geq 1$) is the associated mapping cylender $C_{\pi_i}M_i$
for some fibration $(M_i, \pi _i)$. The metric on $X$ is such that its restriction to each $X_i$ is quasi-siometric to one of the form like $g_1$.

Now given an open Riemannian manifold  $(Y,g)$,  let $A^i =A^i (Y)$ denote the space of $C^{\infty}$ $i$-forms on $Y$, and let $L^2 = L^2 (Y)$ be the completion of $A^i$ with respect to the $L^2$ metric induced by $g$.  Then the exterior differential $d$ has the following domain
\[
\text{dom} (d)=\{\alpha\in A^i (Y)\cap L^2(Y); d\alpha\in L^2 (Y)\}. 
\]
 
The cohomology associated with  the following sequence
\[
...\rightarrow A^{i-1} \xrightarrow{d} A^{i} \xrightarrow{d}
A^{m+1}\rightarrow...
\]

is called $L^2$-cohomology of $Y$, and is denoted by $H^i _{(2)} (Y)$. There is also a natural map $i_{(2)}:H^i _{(2)} (Y) \rightarrow H^i _{(2), \#}(Y)$, which is always an isomorphism between the $L^2$ cohomology of $d$ and  the $L^2$ cohomology of its strong closure 
$\bar{d}$, 
\begin{equation}\label{hdiez}
H^i _{(2), \#} (Y) = \ker \bar{d}_i /\text{im }\bar{d}_{i-1}.
\end{equation}

 Likwise the space of $L^2$ harmonic $i$-forms  is defined by 
\[
\mathcal{H}^i _{(2)} (Y)=\{\theta\in A^i \cap L^2; d\theta=\delta\theta =0\}.
\]
 
Following Cheeger \cite{cheeger1}, strong Hodge theory holds when the Hodge  map
\[
\mathcal{H}^i _{(2)} (Y)\rightarrow H^i _{(2)} (Y),
\]

induces an isomorphism. In this case,  the $L^2$-cohomology of $X$ is of finite dimension, and Stoke's theorem holds on $X$ in the $L^2$ sense, and hence the Hodge theorem holds according to Cheeger and Dai \cite{cheeger3}. 
The topological   ingredient of this cohomology group is also explored through intersection cohomology of  Goresky and MacPherson. 

All these are restricted to  $ L^2$ deRham cohomology for  real manifolds or complex manifolds with standard nonsingular complex structure. We are going to develop a $L^2$  Hodge cohomology version for complex manifolds with singular complex structure.

\subsection{Definition of Degenerate Complex Manifolds}\label{sec42}

%\subsubsection{ Singular complex manifolds}

The following definition is a generalization of the degenerate complex manifolds associated with the  structures   $J'$  in the twistor sphere $\mathbb{S}$  described in section 2. All the theorems in this section, including $L^2$ Hodge decomposition, fundamental inequality and Sobolev embedding theorems  are proved for the  class of singular complex manifolds introduced in the definition (\ref{def2}). Although the weak and strong Hodge theory  discussed in section \ref{weakh} and afterwards are restricted to the above mentioned case.

\begin{defi}\label{def2}
Let $X$ be a smooth real manifold  of dimension 2n and let $D\subset X $  be a real smooth submanifold of real dimension k. Assume that $X\setminus D$ is endowed with a K\"ahler structure $(w',g',J')$ with the following assymptotic behavior near $D$:

1) The K\"ahler metric itself and its Ricci curvature tensor both  are bounded   with respect to  the sup norm  $|.|_{sup}$ induced by the degenerate metric $g'$ in $X\setminus D$
\\

2) Around each point $p\in D$ there exists an open set $U$  and a local coordinates system $\{(x_1,...x_{2n})|x_i\in \mathbb{R},1\leq i\leq 2n\}$ defined on $U$  such that $D\cap U =\pi ^{-1} (0)$  and $det (g') = |r|^s + o(r^s)$. Where $\pi :\mathbb{R}^{2n} \rightarrow \mathbb{R}^{2n-k}$ denotes projection on the  last $2n-k$ coordinates, $r=\sqrt{\sum_{i=2n-k+1} ^{2n} x_i ^2}$, and $g' $ is the singular K\"ahler metric.
\\

3) If $\mathcal{S}$ denotes the sphere bundle of the normal bundle $N$ of $D$ in $X$ then there exists an $\epsilon >0$ and a smooth application 
$T: \mathcal{S}\times [0,\epsilon ) \rightarrow U_{\epsilon} $ onto an open neighborhood $U_{\epsilon}$ of $D$ in $X$   such that $T|_{\mathcal{S}\times ,(0,\epsilon)}$  defines a diffeomorphism onto $U_{\epsilon} \setminus D$,   $T|_{D\times\{0\}} : \mathcal{S}\times \{0\} \rightarrow D $ sends $(x,0)$ to $p(x)$ where $p:N\rightarrow D$ is the bundle projection. And such that $T_* (\frac{\partial}{\partial r}) 
\perp D_{r}$   where for $0<r<\epsilon$,  $D_r =  T (\mathcal{S}\times\{r\})$  and $\frac{\partial}{\partial r}$ is the vector field on $S\times (0,\epsilon)$ tangent to the direction of the interval $(0,\epsilon)$ at each point. 

We also assume that  if a neighborhood of a point $p\in D$ is endowed with a   coordinates system such that  $r$ is one of its coordinates then $\|dr\|_{g'} =r^{t} +o(r^t)$ for some $t>0$.
\\

4) The parameters  $s$  and $t$, defined respectively in  parts  2 and 4, satisfy 
\[
 2n-k+s-2>3t
\]
\\

%5) There exists a neighborhood $U$ of $D$ such that $\partial D$ is strongly pseudoconvex.
\end{defi}

The  assumption (3) about $\|dr\|_{g'}$   implies  that the  total volume of the submanifolds $D _r$  is vanishing at least  of order $O(r^{2n-k+t})$ around $D$  since $\|J' dr \|_{g'}= \|dr\|_{g'}$  and this covector lives in the cotangent bundle of $D _r$.

  It is clear that for the degenerate complex manifolds obtained through degenerate Monge Amp\`ere equation we have $s=2$ and $t=1$ and it  satisfies the conditions of the above  definition.

\subsection{Weighted $L^2$ Spaces and $L^2$ Hodge Decomposition}

Now for a singular K\"ahler manifold   $X$ as in definition (\ref{def2}) we denote by $W'^{p,q}_{l} (X)$ the Sobolev space of $(p,q)$-forms on $X\setminus D$   with respect to the degenerate metric $g'$

\[
 W'^{p,q} _l (X):= \mathcal{H}_l (X\setminus D, \Lambda ^{p,q} (X\setminus D)),
\]

where the sobolev $l$-norm $\|.\|' _{2,l}$  is defined by

\[
\|\phi \|' _{2,l} = \|\phi\|' _2 +\|\nabla_{g'} \phi\|'_2  +...+\|\nabla^l  _{g'}\phi\|'_2.  
\]

Here $\nabla_{g'}$ is the Levi-Civita connection associated to $g'$ on $X\setminus D$, and  $\mathcal{H}_l (M,E)$ denotes the Sobolev space of order $l$ of sections of  a vector bundle $E$ over a manifold $M$.
\\

We represent  by  $L' _2 (X)$,  the
space of functions in $X\setminus D$ which are $L^2$ integrable with respect
to the volume form $d\Omega_{g'}$  induced by the metric $g'$ and  let $L'^{p,q} _2$ (resp. $L'^m
_2$) be the space of $(p,q)$-forms (resp. $m$-forms)  $\alpha$ on $X\setminus D$
s.t. $\int \alpha \wedge
* \alpha < \infty $ .
\\

Also  $W' _1 (X)$  denotes   the space of functions
$f\in L' _2 $  s.t. $df \in L'^1 _2  $.

\begin{lemma}\label{dense}
If $A_0 ^{p,q} (X\setminus D)$ denotes the space of smooth $(p,q)$-forms with compact suppport in $X\setminus D$  then 
$A_0 ^{p,q} (X\setminus D)$ is dense in $W'^{p,q} _l (X)$.
\end{lemma}

\begin{proof} 

Consider a bump function $A:\mathbb{R}\rightarrow \mathbb{R}$ such that

\[
A(r)=\begin{cases}
1 & \text{ if } r\geq b,\\
0&\text{ if } r\leq a,
\end{cases}
\]

where $0<a<b$ are two fixed real numbers.
For $m\in \mathbb{N}$, we  define $a_m : X\rightarrow \mathbb{R}$ as follows

\[
a_m (y,r) =   A(\frac{m}{c^l(t+1)}r^{c(t+1)}).
\]

Here $y\in \mathcal{S}$ and $r>0 $, are as presented in the  definition (\ref{def2}). Also, $c$ is a constant that will be determined later.
%%%%%%%
\\
\\
 
For $\alpha\in A^{p,q} (X\setminus D)$, we have 

\[
\nabla^l _{g'} (a_m \alpha)=\sum b_{i,l} \nabla ^i _{g'} a_m \wedge \nabla^{l-i} _{g'} \alpha,
\]

%according to (\ref{nab1}) and (\ref{nab2})  we have 
%\[
%=\|\nabla ^i _{g'} a_m\|_{g'} = O(m^i)
%\]
 
where  $b_{k,l}$'s are combinatorial constants.  Also for $i\geq 1$, we know that the support of $\nabla ^i _{g'} a_m$ is contained in $U_{\epsilon_1, \epsilon_2}:= U_{\epsilon_2}\setminus U_{\epsilon_1}$, where
    $\epsilon_1 =(\frac{(t+1)a c^l}{m})^{\frac{1}{c(t+1)}}$,   $\epsilon_2 =(\frac{(t+1)bc^l}{m})^{\frac{1}{c(t+1)}}$, and   $U_{\epsilon}$ for $\epsilon >0$ denotes the  neighborhood of $D$ in $X$ given by $r<\epsilon$. We have 
\\

\begin{equation}\label{mbkb}
\|\nabla_{g'} ^l (a_m  \alpha) -\nabla _{g'} ^l \alpha\|_{2} '= \| (a_m -1) \nabla_{g'} ^l \alpha +\sum_{i\geq 1} b_{i,l} \nabla_{g'} ^{i} a_m \wedge  \nabla_{g'} ^{l-i}\alpha\|' _2
\end{equation}

%%%%%%%%%%%

Then the vector field  $v=\frac{1}{r^t}\frac{\partial}{\partial r}$  defined in $U_{\epsilon}$ has bounded $g'$ -norm   of $O(1)$ and we have
\begin{equation}\label{nab1}
D_v a_m (y,r) = \frac{1}{r^t}\frac{\partial}{\partial r} a_m (y,r)= \frac{m  r^{(c-1)(t+1)}}{c^{l-1}}A' (\frac{m}{c^l(t+1)}r^{c(t+1)}).
\end{equation}
 According to  Fa\`a di  Bruno's formula  
\[
(D_v)^l a_m =\sum \mathcal{C}_{\mu_1,...,\mu_l}D_{v} ^{\mu_1+...+\mu_l} A (\frac{m}{c^l(t+1)}r^{c(t+1)})
\Pi _{j=1}  ^l  (D_v ^ j (\frac{m}{c^l(t+1)}r^{c(t+1)})) ^{\mu_j},
\]

where $\mu_,...,\mu_l$  are nonnegative integers satisfying $\sum j\mu_j =l$,  and $\mathcal{C}_{\mu_1,...,\mu_l}$ are constants.
So

\[
 (D_v ^ j (\frac{m}{c^l(t+1)}r^{c(t+1)})) ^{\mu_j}=\bigg (\frac{m}{c^{l}} (t+1)^{j-1}\left [\prod_{i=0} ^{j-1} (c-i)\right ] r^{(c-j)(t+1)}\bigg )^{\mu_j}\leq\kappa m^{\frac{j\mu_j}{c}}\leq  \kappa m^{\frac{l}{c}},
\]

for a  constant $\kappa$ which depends on $a, b$, $t$ and $c$.
Here in the first  inequality we are using the definition of $\epsilon_1$ and $\epsilon_2$ and the last inequality comes from the above
properties of $\mu_1,...\mu_l$.   We also have

\[
Vol (U_{\epsilon_2} \setminus U_{\epsilon_1}) =O((\frac{1}{m})^{s+2n-k })
\]
and

\begin{equation}\label{nab2}
\frac{1}{r^{t_i}}\frac{\partial}{\partial y_i}a_m (y,r)=0,
\end{equation}

where $\frac{\partial}{\partial y_i}$ is a vector tangent to $\mathcal{S}$ and $t_i$ is such that, $\|\frac{\partial}{\partial y_i}\|_{g'}=O(r^{t_i})$.
%%%%%%%

Thus we obtain

\[
\|\nabla _{g'} ^i a_m \wedge \nabla^{l-i} _{g'} \alpha\|_2 '  \leq \kappa ' m^{\frac{l}{c}-s-2n+k},
\]

where $\kappa'$ is a constant which only depends on $a,b,t$ and $c$.

Due to property (4) of definition (\ref{def2}), this means that if $c>l$  and $i\geq 1$, we  will have
\[
\lim_{m\rightarrow \infty}  \|\nabla _{g'} ^i a_m \wedge \nabla^{l-i} _{g'} \alpha\|_2 ' =0.
\]
It can also be seen that 
\[
\lim_{\epsilon\rightarrow\infty } \| (a_m -1) \nabla_{g'} ^l \alpha \|_2 '=0.
\]

 Hence  due to (\ref{mbkb}) we get

\[
\lim_{m\rightarrow \infty}\|\nabla_{g'} ^l (a_m  \alpha) -\nabla _{g'} ^l \alpha\|_{2} '=0.
\]
 \end{proof}
\begin{lemma}\label{k}
If $f\in L'^{2} (X)$ is such that $lim_{r\rightarrow 0} \int_{  D_r} fd\Omega' (r) $ exists, then this limit must be zero: 
\[
lim_{r\rightarrow 0} \int_{D _r}  fd\Omega' (r)=0
\] 
Here $d\Omega' (r)$  denotes the volume form induced  by the metric $g'$ on $D _r$  .
\end{lemma}

\emph{Proof:}  Let $D$, $D_r$ and $\mathcal{S}$ be as in definition (\ref{def2}). Let $d\Omega_N $ denote an ordinary volume form on $\mathcal{S}$ iduced by a smooth riemaniann metric.  Then we have $d\Omega' (r)= vr^{2n-k-1}d\Omega _N$, where   $v$ is such that  $vr^t =O(r^s)$, and we know that  
\begin{equation}\label{a}
 \int f^2r^t dr d\Omega'(r)= \int f^2 r^{t+2n-k-1}v dr d\Omega_N <\infty,
\end{equation}

and  that  the following limit exists
\begin{equation}\label{b}
lim _{r\rightarrow 0} \int f d\Omega' (r) = lim _{r\rightarrow 0} \int fvr^{2n-k-1}d\Omega_N.
\end{equation}

Thus by Cauchy-Schwarz inequality  we obtain

%\begin{equation}\label{eq100}
%|\int fd\Omega ' (r) |^2=|\int fvr^{2n-k-1} d\Omega _N |^2 \leq c \int f^2  v^2 r^{2(2n-k-1)} d\Omega_N ,
%\end{equation}

%for some constant $c$ which depends on $d\Omega_N$.
%One can also deduce that
%\[
%\begin{matrix}
%\int |\int fd\Omega '(r) |^2 r^t dr \leq c \int \int f^2 v^2 r^{2(2n-k-1)+t} d\Omega_N  dr &\leq& c' \int f^2vr^{2(2n-k-1)+s} d\Omega_N dr 
%\\
%&=&c'\int( f^2 v r^{t+2n-k-1})  r^{s-t+2n-k-1} d\Omega_N dr

%\end{matrix}
%\]

\[
\begin{split}
\int |\int fd\Omega '(r) | r^t dr &=  \int| \int f v r^{(2n-k-1)} d\Omega_N| r^t  dr 
\\
\quad &=  \int| \int f \sqrt{v} r^{\frac{(2n-k-1)}{2}} \sqrt{v} r^{\frac{(2n-k-1)}{2}} d\Omega_N| r^t  dr \\
\quad &\leq  \tilde{c} \int| \int f \sqrt{v} r^{\frac{(2n+t-k-1)}{2}}   d\Omega_N| r^{\frac{(2n-k-1)}{2}+\frac{s-t}{2}}  dr\\
\quad &\leq \tilde{c}\left ( \int |\int  f \sqrt{v} r^{\frac{(2n+t-k-1)}{2}}   d\Omega_N |^2  dr\times \int  r^{(2n-k-1)+s-t}  dr  \right ) ^\frac{1}{2}\\
\quad & \leq \left (\tilde{c} c'\int |\int  f^2 v r^{(2n+t-k-1)}  dr d\Omega_N | \times  \int  r^{(2n-k-1)+s-t}  dr \right ) ^{1/2}\\
\quad & =O(r^{\frac{2n-k+s-t}{2}}).
\end{split}
\]
Here in the third line  we are using the relation $vr^t =O(r^s)$ from which it follows that $vr^t\leq \tilde{c}r^s$ for some constant $\tilde{c}$. Also in the forht and fifth line we apply Cauchy-Schwarz inequality, where $c'=\int d\Omega_N$ and for the last line we apply (\ref{a}).  Now if the  value of the limit in relation (\ref{b})  is different from zero then  by comparing  the order of zeros in terms of $r$  in the two sides of the above relation     we reach to a contradiction 
 if $ 2n-k+s-t>2t+2$ or equivalently  if $ 2n-k+s-2>3t$.  And this holds according to part (4) of definition (\ref{def2}).
\\
\\
%For a function $f\in   W' _1$ in local coordinates around a
%point $p\in D$ the derivatives of $f$ with respect to
%$\partial _u,
%\partial _v, \frac{1}{r}\partial _r $ and $\frac{1}{r ^2}
%\partial _{\phi}$ belong  to $L' _{2, {loc}}$, where $(r,\theta)$ is the polar coordinates on $(x,y)$ plane
%passing through $p$.
%
 %It can be easily verified that if $X$ is a degenerate K\"ahler manifold in the sense of definition (\ref{def1}) then  the Sobolev %space $W' _1(X)$ coincides with the
%ordinary Sobolev space $H_1 (X)$.\\
  \\
\\

\begin{cor}\label{ck}

i) If $\eta\in W'^{2n-1} _1$, then $\int_X d\eta =0$.

\hspace{2cm}ii )If $\eta\in A^{2n-1} (X\setminus D) \cap L'^{2n-1} _2$ is such that $d\eta\in L'^2$,  then $\int_X d\eta =0$

\hspace{2cm}iii) If $\alpha\in W'^p _1$, $\beta\in W'^q _1$  and $p+q=2n-1$, then

\[
\int_X d\alpha\wedge \beta = (-1)^{p} \int_X \alpha\wedge d\beta 
\]
 
\end{cor}

\begin{proof}
For a proof  of (i) it suffices to define $f$ in a neighborhood $U$ of $D$ by $j_r ^* \eta =fd\Omega'  (r)$. Then it can be seen that $f\in L'_2 (U)$. By Stokes theorem the limit $\lim _{r\rightarrow 0}\int_{D_r} fd\Omega'(r) $ exists. Then one can apply the above lemma (\ref{k}). The proof of the other parts is similar.
\end{proof}

 Clearly $W'^{p,q} _{1} \subset L'^{p,q} _2 =W'^{p,q} _0$, and  the following
$L^2$-decomposition holds:
\begin{theo}\label{l2dec}
\[
L'^{p,q} _2 = \mathcal{H} _2 ^{p,q} \oplus
\overline{\bar\partial _{J'} (A^{p,q} (X\setminus D) \cap L'^{p,q} _2  )}\oplus \overline
{\bar\partial _{J'} ^*  (A^{p,q}(X\setminus D)\cap L'^{p,q}_2 )} =V_0 \oplus V_1 \oplus V_2,
\]

where the space $\mathcal{H} _2 ^{p,q}$ is defined by $
\mathcal{H} _2 ^{p,q} := \{\alpha \in L'^{p,q} _2 | \
\bar\partial_{J'} \alpha = \bar\partial ^* _{J'} \alpha = 0 \}$.
 \end{theo}

%Using the inequality $\|. \| _ 2 \leq c \|.\| _{2,1}$  we can see that $W_i$, $i=1,2$ is a closed subspace of $W' ^{p,q} _{1}$ %and we have
%
%\[
%W'^{p,q} _{1} = {\mathcal{H} ^{p,q} \oplus W_1 \oplus
%W_2 }
%\]
%is a decomposition to closed subspaces with respect to $\|.\| _{2,1}$-norm.

According to lemma (\ref{dense}), we also have

\begin{lemma}
\[
V_1  =\overline{ \bar\partial _{J'} A^{p,q} _0 (X\setminus D)}, \hspace{0.5cm}
V_2 =\overline{ \bar\partial _{J'} ^* A^{p,q} _0 (X\setminus D)},
%
%V_i \cap  W^{p,q} _{1,2} =\overline{W^{p,q}_{1,2} \cap \bar\partial _{J'} A^{p,q}}
\]

where $A^{p,q} _0 (X\setminus D)$ denotes the space of smooth  $(p,q)$ forms with compact support in $X\setminus D$.
\end{lemma}

 \subsection{Fundamental Inequality.} 
In this subsection we would like to prove  fundamental inequality for the Laplace operator $\Delta_{g'}$ of a singular complex manifold in the sense of definition \ref{def2}.
%We start by the following lemma:
%\begin{lemma}
%The space $A^{k} _0 (X\setminus D)$ is dense in the Sobolev space $W^k _2 (X\setminusD)$.
%\end{lemma}
%
%\emph{Proof:} Take an element $\alpha\in W^k _2$, and let $\epsilon$ be small such that $|\alpha|_{U_{\epsilon}}  |' _2$ be %small. Let $\psi\in A^k (X\setminus U_{2\epsilon})$ be such that $|\alpha-\psi|' _{2}$ be small.  Consider a bump function %$a_{\epsilon}$ such that 
%$a_{\epsilon }|_{X\setminus{U_{2\epsilon}}} \equiv 1$ and $a_{\epsilon }|_{U_{\epsilon}}\equiv 0$ (constructed as in the %proof of previous lemma). We claim that $\psi =a_{\epsilon} \alpha$ is close to $\alpha$ in $W' ^k _2$. To see this all we %need is that $a_{\epsilon }$ can be chosen such that $|a_{\epsilon}|_{U_{2\epsilon}}|' _2$ becomes small.\square
%\\

\begin{theo}
On a singular complex manifold in the sense of definition (\ref{def2}) and for $\psi\in W'^{p,q} _2$ there exists a constant $c$ such that ,

\[
\|\Delta_{g'}\psi\|' _2\leq c\|\psi \|'_{2}. 
\]
\end{theo}
\begin{proof}
 Suppose that $\phi_1,...,\phi_n$ is a local unitary coframe and consider a form $\psi\in A^{p,q} _0 (X\setminus D)$ which is locally written in the form 
\[
\psi =\frac{1}{p!}\frac{1}{q!} \sum _{I,J} \psi_{I,J} \phi _I \wedge \bar{\phi}_J.
\]
We also set $\Phi ' =\phi_1\wedge...\wedge \phi_n$. Now following standard computation %\cite{GH} page 99,
 modulo first order terms we have

\begin{equation}\label{W2}
2\Delta_{g'} \psi \wedge * \Delta_{g'} \psi= d\eta _1 - d\eta _2 +d\eta_3 + 2\sum |\psi_{IJij}|^2 w^n  +A^1 (\psi ),
\end{equation}
 where 
 
\[
\eta_1 = \sum (-1) ^i \psi_{IJ \bar{i}}\bar{\psi}_{IJ j\bar{j}} \phi_1 \wedge ... \hat{\phi_i}\wedge ... \phi_n \wedge \Phi ' =<\bar{\nabla} \psi, \Delta_{g'} \psi> \wedge \omega '^{n-1}
\]
\[
\eta_2 = \sum (\psi_{IJ\bar{i}}\bar{\psi}_{IJ \bar{i}j} +  \psi_{IJ\bar{i}\bar{j}}\bar{\psi}_{IJ\bar{i}})  \phi_1 \wedge ... \hat{\phi_j}\wedge ... \phi_n \wedge \Phi '=\sum_i  \bar{\nabla} |\psi_{IJ i}|^2 \wedge \omega '^{n-1}
\]

\[
\eta _3 =\sum ( \psi_{IJ\bar{j}j}\bar{\psi}_{IJ\bar{i}})   \bar{\psi_1}\wedge ...\hat{\bar{\psi}_i}\wedge ... \bar{\psi}_n \wedge \psi ' = <\Delta_{g'} \psi, \nabla \psi> \wedge \omega '^{n-1},
\]
  
and where in the K\"ahler case,  $A^1(\psi)$ is an algebraic operator depending on the Ricci curvature. 
On the other hand $\eta = \eta_1 -\eta_2 +\eta_3$ has  its support in $X\setminus  D$ and thus we have $\int_{X} d\eta =0$.

Since according to  definition \ref{def2}, the ricci curvature is assumed to be bounded then by  using  lemma (\ref{dense})     we obtain te lemma.

\end{proof}

One can then  conclude  that
 
\begin{cor}
The operator $\Delta_{g'} :W'^{p,q} _{2} \rightarrow W'^{p,q} _{0}$ is a bounded Fredholm operator.

\end{cor}

\begin{theo}\label{reli}
The Sobolev embedding $j: W'^{p,q} _1 \hookrightarrow L'^{p,q} _2$  is  compact.
\end{theo}

\begin{proof} Consider a bounded sequence $\alpha_n$ of the elements of $W'^{p,q} _1$. First we claim that for each $\delta >0$ there exists an $\epsilon >0$ such that  for all $n\in \mathbb{N}$, we have $\|\alpha_n|_{U_{\epsilon}} \|' _2< \delta$. To see this, we consider  for each  $0<s<r $ and for $x\in \partial U _s$,   a  differential form  $\eta _n ^s$  defined by $\eta _n ^s (x) := \|\alpha _n (x)\|_{g'} ^2 d\Omega ^s _{g'}(x)   $  where $d\Omega^s _{g'} (.)$ is the volume form on $\partial U_s$ induced by the metric $g'$. One can see that
\[
|\int _{\partial U_{r}} \eta_n| = |\int _{U_{r}} d\eta_n|\leq c,
\]

where the equality comes from part (ii) of corollary  (\ref{ck}), and the inequality is a consequence of the existence of a uniform bound on $\|\alpha_n\|'_{2,1}$ for all $n\in \mathbb{N}$. ($\|\|'_{2,*}$ is the norm of $W'^{p,q} _{*}$).
Thus $c$ is a constant which   is independent of $n$ and $r$.
%Also $U_r$ is a neighborhood of $D$ as in definition (\ref{def2}).
\\
Hence we have 
\[
(\|\alpha_n|_{U_R}\| ' _2)^2= \int_{0\leq r\leq R}  (\int _{\partial U_r} |\alpha_n|^2 d\Omega^r _{g'}) r^t dr   \leq c\frac{R^{t+1}}{t+1}.
\]
 So if $c\frac{\delta ^{t+1}}{t+1} <\epsilon ^2$  then the inequality  $\| \alpha_n |_{U_{\delta}}\|' <\epsilon$  holds for all $n\in \mathbb{N}$. In order  to complete the proof of lemma, we consider  a sequence of neighborhoods $\{ U_{\delta_m}\}$ of $D$ such that 
\begin{equation}\label{den}
\|\alpha_n| _{U_{\delta_m}}\|' _2 <\frac{1}{m}, \hspace{0.5cm} \text{ for all } n\in \mathbb{N}.
\end{equation}

By applying standard Rellich compactness theorem on $X\setminus U_{\delta_m}$,  we can find a subsequence $\{\alpha_{m_n}\}$ converging in $L'^2$ to an element $\beta_m \in L'^{p,q} _2 (X\setminus U_{\delta_m})$. Thus for the diagonal subsequence $\{\alpha_{n_n} \}$, we know that 
\[
\alpha_{n_n} |_{X\setminus U_{\delta_n}} \rightarrow \beta,
\]
for an elemen  $\beta\in \cap_m  L'^{p,q} _2 (X\setminus D)$   which is defined by $\beta|_{X\setminus U_{\delta_m}} = \beta_m$. On the other hand 
by (\ref{den})
\[
\|\beta |_{U_{\delta_{m_1}} \setminus U_{\delta _{m_2}}}\|'  _2 \leq lim_{n\rightarrow \infty} |\alpha_{n_n} |_{U_{\delta_{m_1}} \setminus U_{\delta _{m_2}}} |' \leq \frac{1}{m_1},
\]

for all $m_1<m_2 \in \mathbb{N}$. So by taking the limit $m_2\rightarrow \infty$ we obtain $\|\beta | _{U_{\delta _{m_2}}} \|_2 ' <\frac{1}{m_1}$. This means that $\beta\in L'^{p,q} _2 (X)$  and the proof is completed.
 \end{proof}

\section{ Weak Hodge Theorem}
\label{weakh}
In order to derive weak and strong Hodge theory we assume that $(X,J')$ is the singular complex manifold described in  section \ref{sec2} which is constructed through the solutions to DCMA equations.

Let us define

\begin{equation}\label{C}
C^{p,q} (X):= A^{p,q} (X\setminus D ) \cap L'^{p,q} _2 (X),
\end{equation}
%be the space of $(p,q)$-forms with singularities  along $D$.  
%We also need the following definitions:
%\begin{defi}
%The space $\mathcal{H}^{p,q} (X)$ is defined as 
%\[
%\mathcal{H}^{p,q}(X)=\{\alpha\in\mathcal{H}^{p,q} _2|\alpha +\bar{\partial} _{J'}\beta \in W'^{p,q} _1 (X) \text{ for some } %\]
%\end{defi}

and  consider    the following $\bar\partial _{J'}$ -complex of differential forms:
 
\begin{equation}\label{cpq}
...\rightarrow C^{p,q-1} \xrightarrow{\bar\partial_{J'}} C^{p,q}
\xrightarrow{\bar\partial_{J'}} C^{p,q+1}\rightarrow...
\end{equation}
 
Obviously the domain of definition of $\bar\partial _{J'}$ in the above complex is strictly smaller than
$C^{*,*}$.

Then we have the following theorem:

\begin{theo}\label{coh2}
For every $\bar\partial _{J'}$-closed $(p,q)$-form $\alpha\in C^{p,q} (X)$,  there exists a unique harmonique form $\alpha _0 \in \mathcal{H}^{p,q} _2$  and $\beta\in C^{p,q-1} (X )$ such that, $\alpha=\alpha_0 + \bar\partial _{J'}\beta$. Moreover if $\alpha\in W'^{p,q} _1$, then $\beta$ can be chosen to belong to $W'^{p,q-1} _1$ as well.
 
\end{theo}

 \begin{proof}
  Assume that  $\alpha$ is a
$\bar\partial _{J'}$-closed $(p,q)$- form such that $\alpha \in W'^{p,q} _{1}$  (or if $\bar{\partial}_{J'} ^* \alpha \in L'^2$). Then by theorem (\ref{l2dec}),  one can find an
element $\alpha _0 \in \mathcal{H} ^{p,q} _2$ and a sequence
$\{\theta _n \}_n$ of the elements of $A_0 ^{p,q}(X\setminus D)$ such that
\[
 \alpha_0
+\bar\partial _{J'} \theta _n \rightarrow \alpha \text{ in } L'_2.
\]
  If we consider a decomposition of $\theta_n$ of the form

\[
\theta_n =\theta_n ^h + \theta_n ^1 +\theta_n^2,\hspace{0.5cm} \theta_n ^i \in V_i \hspace{0.5cm } \text{for } i=0,1,2
\] 

 ($V_i$'s are defined in theorem (\ref{l2dec})),   then we can find sequences  $\eta_n ^m , \beta_n ^m \in A_0 ^{p,q} (X\setminus D )$  such that
\[
 lim _{m\rightarrow {\infty}}\bar\partial_{J'}\eta_{n} ^m  \rightarrow \theta _n ^1, \hspace{0.5cm} 
 \text{and } \hspace{0.5cm}\ lim _{m\rightarrow {\infty}}\bar{\partial_{J'} }^* \beta_n ^m \rightarrow \theta_n^2
\] 

 in $L'^2$. Thus $\alpha_0 +\bar{\partial}_{J'}\bar\partial_{J'} ^* \beta_n ^m$ converges  in the weak sense  towards $\alpha_0 +\bar{\partial} _{J'} \theta_n$ as $m$ grows to infinity . 

Replacing $\theta_n$ by some appropriate subsequence  $\bar{\partial} ^* _{J'}\beta^{m_n} _n$,  we can  assume   that $\bar\partial _{J'} ^* \theta _n =0$.  Hence one can  conclude that
\begin{equation}\label{dthe}
\Delta_{g'} \theta _n = \bar\partial _{J'} ^* \bar\partial _{J'}
\theta _n \rightarrow \bar \partial _{J'} ^* ( \alpha-\alpha _0 )=\bar\partial ^* _{J'} \alpha,
\end{equation}
in the weak sense (in $L'^{p,q} _2$).   It is worthy to emphasize  that here we are using the fact that $\alpha\in W'^{p,q} _1$.(or $\bar{\partial}_{J'} ^* \alpha \in L'^2$). We have 
% since $\theta_n \in A^{p,q} (X\setminus D)\cap L'^{p,q}_2$ we can assume that $\Delta_{g'} \theta_n$ tends to $\bar %\partial _{J'} ^* ( \alpha-\alpha _0 )$ in $L'^{p,q} _2$ on $X\setminus U_{\frac{1}{m}}$. 

%Since right hand  side belongs to $\L'^{p,q}_2$ we can find a subsequence such that $\Delta_{g'}\theta_n \rightarrow \bar |%\partial _{J'} ^* ( \alpha-\alpha _0 )$ in $L'^2$ on $X\setminus D$.

\begin{equation}\label{preponc}
\|\nabla_{g'}\theta_n\|'^2 _{2}=
|\int \langle  \Delta_{g'} \theta_n , \theta_n \rangle_{g'}| =|\int \langle \bar\partial^*_{J'} \alpha, \theta_n \rangle_{g'}|\leq  \| \bar\partial^*_{J'} \alpha \|_{L'^2} \|\theta_n \|_{L'^2},
\end{equation}

where the first equality is proved through integration by part, the second equality is deduced from (\ref{dthe}), and the last inequality results from the Cauchy-Schwartz inequality.

Now let $\{U_i\}_{1\leq i\leq K}$ be a  finite  open covering of $X$ such that $U_1$ is a neighborhood of $D$ which is  Stein, in the sense of section (\ref{solv}), and $U_2,...,U_K$  have positive distance from $D$ and are domains of holomorphy for example holomorphic discs in
$\mathbb{C}^2$.
From the inequality (\ref{3990}) we have
\begin{equation}\label{hor1} 
\|\theta_n|_{U_i}\|'^2 _{\chi (\phi_i)}\leq  (\|\bar{\partial}^*_{J'} \theta_n |_{U_i} \|'^2 _{\chi (\phi_i)} +(\|\bar{\partial}_{J'} \theta_n|_{U_i} \| '^2 _{\chi (\phi_i)}.
\end{equation}

Here $\phi_i$'s are appropriate plurisubharmonic functions on $U_i$'s and $\chi$ is an appropriate convex increasing map.

Since $\chi (\phi_i)$ tends to zero near the boundary, we can take a finer coverig $\{V_i\}_{1\leq i\leq K}$ with $\bar{V}_i\subset U_i$
 and such that
\begin{equation}\label{pl}
\|\theta_n|_{V_i} \|'^2_2 \leq C'\|\theta_n|_{U_i}\|'^2 _{\chi (\phi_i)}
\end{equation}
for an appropriate constant $C'$.

By Cauchy-Schwartz we have 
\begin{equation}\label{p1}
\sum \|\theta_n|_{V_i}\|'^2_2 \geq \frac{1}{K} (\sum  \|\theta_n|_{V_i}\|_2' )^2\geq \frac{1}{K}\|\theta_n \|'^2 _2,
\end{equation}
and from (\ref{preponc}), (\ref{hor1}),(\ref{pl}) and (\ref{p1}) it follows that
\[
\| \theta _n\|'  _2\leq C,
\]
or some constant $C$ independent of $n$. 
 This proves that the sequaent $\{\theta_n\}$ is a uniformly bounded  
 sequence with respect to the induced norm from $L'^{p,q} _2$.
 
One can see  that the sequence $\{\theta_n \}_{n}$ can be taken to be uniformly bounded in $L'_2$, and it follows from  (\ref{preponc}) that  $\{\theta_n\} _n$ is uniformly bounded in $W' _{1,2}$. Thus by Rellich theorem (see theorem  (\ref{reli})) it contains a  subsequence  $\{\theta_{n_j}\}_{j} $, convergent   in $L' _2$ to some $\theta$. Also according to (\ref{preponc}),  $\nabla_{g'} \theta_{n_j}$ is  a Cauchy sequence and hence it is convergent in $L'^2$ as well. This means that $\theta\in W'^{p.q}_1$.  
% If $\theta$ denotes the limit of this subsequence in $L' _2$ then the identity
% $\Delta _{g'} \theta = \bar \partial _{J'} ^* ( \alpha-\alpha _0 )$ holds in the weak sense. So again from (\ref{preponc}) one  can conclude  that $\theta\in W' ^{p,q} _1$.   

Also by standard (non-degenerate) regularity arguments  we can prove the smoothness of $\theta$ outside $D$
in the case where $\alpha\in A^{p,q} (X\setminus D )\cap L'^{p,q} _2$.
%r
%According to orem 2   the operator $\Delta _{g'}$ is bounded and this  implies that   the sequence $\theta _n$
%s also bounded. Thus  we obtain a subsequence converging in the weak sense to an element $\theta \in W'^ {p,q}
%_2$. This is due to the compactness of $W'^{p,q} _2$ in $ H^{p,q} _{(2,0,1)}$. By standard (non-degenerate) regularity 
%arguments  we can prove the smoothness of $\theta$ outside $D$
%in the case where $\alpha\in \A^{p,q} (X)$, although nothing can be said about the regularity along $D$.
%
%
%
%Since $\alpha _1 \in W^{p,q} _{1,2}$ so  $\rho ^2 \bar\partial ^*
%\alpha _1 \in L^2 $ (where $L^2$ is the space constructed via the
%ordinary (nonsingular ) Lebesgue measure. Hence by the theorem ??
%on boundedness of the operator $P$, the sequence $\{\theta _n \}
%_n$ is uniformly bounded in $W ^{p,q-1} _{1,2}$ and thus contains
%a subsequence $\{\theta _ {n_i}\}$which converges in the weak
%sense to some $\theta \in W ^{p,q-1} _{1,2}$.
%
%Now we claim that $\theta $ is smooth outside $D$ but
%nothing can be said about the regularity along $D$. This is
%because when calculating the application of $P $ on $\theta$ we
%have

Consequently $\theta_n\rightarrow \theta$ and $\bar{\partial}_{J'} \theta_n\rightarrow \alpha-\alpha_0$ both  in $L'_2$, from which it follows that 

\[
\langle \theta_n, \bar{\partial} _{J'} ^* \phi \rangle \rightarrow \langle \alpha-\alpha_0 , \phi\rangle  \hspace{0.5cm} \text{ and } \hspace{0.5cm}
  \langle \theta_n, \bar{\partial }_{J'} ^* \phi \rangle \rightarrow \langle \theta, \bar{\partial}^* _{J'} \phi\rangle 
\]
  
for any test function $\phi$. Therefor we get

\[
\bar{\partial}_{J'} \theta=\alpha-\alpha_0. 
\]

  This completes the proof of the second part of the theorem.

% If $\alpha\in C^{p,q} (X)$ belongs also to $W'^{p,q} _1$ the proof is complete by our previous discussion.

In order to prove the first part, by the solvability of  $\bar{\partial}_{J'}$-equation in  a neighborhood of $D$ we can find $\gamma\in C^{p,q-1}$
such that $\alpha-\bar{\partial}_{J'} \gamma \in W'^{p,q} _1$ and we can repeat our previous argument.

\end{proof}

\begin{cor}
If $\alpha\in W'^{p,q} _1$ then in the Hodge decomposition $\alpha=\alpha^h+\alpha^1+\alpha^2$ where $\alpha^h\in\mathcal{H}_2 ^{p,q}$ and $\alpha^i\in V_i$, for $i=1,2$, we have in fact $\alpha^1\in W'^{p,q-1} _1$ and $\alpha^1\in W'^{p,q+1} _1$
\end{cor}
\begin{proof}
As in the proof in section 4.4.1, we find a sequence $\alpha_0+\bar{\partial}_{J'} \theta_n +\bar{\partial}^*_{J'} \eta_n$ such that 
$\bar{\partial}^*_{J'} \theta_n =\bar{\partial}^*_{J'}\eta_n =0$ and which converges in $L'_2$ to $\alpha$. Then $\Delta_{g'} \theta_n\rightarrow \bar{\partial} ^* _{J'} \alpha$ and $\Delta_{g'} \eta_n \rightarrow \bar{\partial}_{J'} \alpha$ in the weak sense in $L'^{p,q} _2$. The rest of argument goes as before.
\end{proof}

The cohomology group associated with the  complex (\ref{cpq}) is denoted by $H^{p,q} _{(2)} (X)$. From the above theorem it follows that
$H^{p,q} _{(2)} (X)\simeq \mathcal{H}^{p,q} _2$.
%We can also consider simila $d$ complex and the associated dRham cohomolgy is denoted by $H^m _{(2)} (X)$. The study of $H^m$ %is th esubject of sectio 4 on Hodge structures.

As in (\cite{cheeger3})  given by the definition (\ref{hdiez},) we can also define  $H^{p,q} _{(2), \#} $  the cohomology associated with the closure $\bar{\partial} ^c _{J'}$ of the operator $\bar{\partial}_{J'}$ . From theorem (\ref{l2dec}) and theorem (\ref{coh2}) it follows that

\begin{cor}
We have the isomorphism
\[
H^{p,q} _{(2)} (X) \simeq H^{p,q} _{(2), \#}  (X)
\]

where  $H^{p,q} _{(2), \#} $ denotes the cohomology associated with the closure $\bar{\partial} ^c _{J'}$ of the operator $\bar{\partial}_{J'}$ 
\end{cor}

%The proof of
%the following theorem is a result of the above discussion:
%
%\begin{theo} (Strong Hodge Theorem)
%For each $\bar\partial_{J'}$-closed $(p,q)$-form $\alpha \in
%C^{p,q} (X) $ there exists  a unique  $\alpha_0 \in \mathcal{H}
%^{p,q}$ and $\beta \in C^{p,q-1} (X)$ s.t. $\alpha =\alpha_0
%-\bar\partial_{J'} \beta$. Thus $H_{C} ^{p,q} (X) \simeq
%\mathcal{H} ^{p,q}$.
%\end{theo}

\section{ Strong Hodge Theory for Degenerate Complex Manifolds}\label{sec6}

To each $\alpha\in A^{p,q} (X\setminus D)\cap L'^{p,q} _1$ we associate a linear map $T_{\alpha}: A^{n-p,n-q} (X)\rightarrow \mathbb{C}$ defined by
\[
T_{\alpha } (\psi ) =\int_X \alpha\wedge \psi.
\]
 
 1) Let $\alpha\in A^{p,q} (X\setminus D)\cap L'^{p,q} _2$ be such that $\bar{\partial}_{J'} \alpha \in A^{p+1,q} (X\setminus D)\cap L'^{p+1,q} _2$.  Here  $\bar{\partial}_{J'} \alpha$ is defined  on $X\setminus D$ by taking ordinary derivatives. Then we have

\begin{equation}\label{deb}
\bar{\partial}_{J'}T_{\alpha} =T_{\bar{\partial}_{J'}\alpha}.
\end{equation}

This is because if $\phi\in A^{n-p ,n-q-1}(X)$ then 
\[
\int_X \alpha \wedge \bar{\partial}_{J'} \phi = (-1)^{p+q}\int_X \bar{\partial}_{J'}\alpha \wedge \phi+
\lim_{\epsilon\rightarrow 0}\int_{\partial U_{\epsilon}}    j^* (\alpha\wedge \phi ).
\]

In fact in the above relation the first two integrals exists because $\alpha\in  L'^{p,q} _2$ and  $\bar{\partial}_{J'} \alpha \in L'^{p+1,q} _2$. Hence the last limit  exists as well. We can now apply  lemma  \ref{k}  to conclude that this limit  should vanish.
 Therefore   (\ref{deb}) is proved.
\\

2) We can deduce from the above observation that
if $\alpha\in A^{p,q} (X\setminus D)\cap L'^{p,q} _2$ satisfies

\begin{equation}\label{eq21}
 \bar{\partial}_{J'}\alpha\in A^{p,q+1} (X\setminus D)\cap L'^{p,q+1} _2,  \hspace{0.5cm} \bar{\partial}_{J'} ^* \alpha\in A^{p-1,q} (X\setminus D)\cap L'^{p-1,q} _2,
\end{equation}

and 

\begin{equation}\label{eq22}
\bar{\partial}_{J'} \bar{\partial}_{J'}^* \alpha , \bar{\partial}_{J'}^* \bar{\partial}_{J'}  \alpha\in  A^{p,q} (X\setminus D)\cap L'^{p,q} _2,
\end{equation}

(here all the derivatives $ \bar{\partial}_{J'}$ and $\bar{\partial}_{J'}^* $ are ordinary derivatives on $X\setminus D$.)
Then we have
\begin{equation}\label{deltaa}
\int_X \Delta_{g'} (\alpha) \wedge\phi = \int_X \alpha\wedge \Delta_{g'} (\phi),  
\end{equation}

or in other words
\[
T_{\Delta_{g'} (\alpha) }=\Delta_{g'}  T_{\alpha}.
\]

In particular if 
\[
 \bar{\partial}_{J'} \alpha = \bar{\partial}_{J'} ^* \alpha =0,
\]

then since all the equations \ref{eq21} and \ref{eq22} hold  true then in the weak sense with respect to the degenerate volume element $d\Omega_{g'}$ on $X$, we obtain
\[
\Delta_{g'} (\alpha)  =0
\]
% associated with the degenerate metric $g'$.
\\

3) Conversely, assume that $\Delta_{g'} \alpha =0$  for some  $\alpha \in A^{p,q} (X\setminus D)\cap L'^{p,q} _2 $  in the weak  sense as of the  relation (\ref{deltaa}). Then it is not difficult to see that
 the operator $|\sigma|^{2}\Delta_{g'}$ has the type of the operators  of Grushin  ( in \cite{gru3}  proposition 1.10),  and  this relation is equivalent to say that 
\[
|S|^{2}\Delta_{g'}\alpha =0
\]

in the weak sense with respect to the ordinary measure $\frac{1}{|S|^s} \mu_{g'}$.  It then follows that $\alpha$  must be smooth  on $X$.  Here $S$ is  a holomorphic section  of $[D]$ vanishing along $D$. 

One
 can then deduce that
\[
\langle \Delta_{g'}\alpha , \alpha \rangle =\int_X \Delta_{g'}\alpha\wedge *\alpha = \int \bar{\partial}_{J'} \alpha \wedge * \bar{\partial}_{J'} \alpha+ \int \bar{\partial}_{J'}  ^* \alpha \wedge * \bar{\partial}_{J'} ^* \alpha.
\]

Here we are using (\ref{deb}). We can thus prove the following theorem

\begin{theo}
For the degenerate complex manifold  $(X,J')$, if $\alpha\in  A^{p,q} (X\setminus D)\cap L'^{p,q} _2$ satisfies $\Delta_{g'}(\alpha) =0$ in the weak sense w.r.t. the measure $\mu_{g'}$, then $\alpha$ is smooth $\alpha\in A^{p,q} (X)$, and 
\begin{equation}
\bar{\partial}_{J'} \alpha=\bar{\partial}_{J'} ^*\alpha=0.
\end{equation}
\end{theo}

Since the same argument holds for   $\bar{\partial}_{J'} $ and $\bar{\partial}_{J'} ^*$ replaced by $d$ and $d^*$ respectively
we can conclude that
\begin{cor}
For the degenerate complex manifold $(X,J')$ we  have
\[
\begin{split}
\{\alpha\in&  \oplus_{p+q=m} A^{p,q} (X\setminus D)\cap L'^{p,q} _2| \bar{\partial}_{J'}\alpha =\bar{\partial}_{J'}^* \alpha =0\}
=\{\alpha\in  A^{m} (X\setminus D)\cap L'^{m} _2| d\alpha =d^* \alpha =0\}\\
\quad &= \{  \alpha\in  A^{m} (X\setminus D)\cap L'^{m} _2| \Delta_{g'} \alpha =0\}. 
\end{split}
\]
\end{cor}

\section{Singular $\partial\bar{\partial}$- Lemma}\label{secb}

\begin{lemma}\label{ddbr}($\partial\bar{\partial}$)
Let $\alpha\in C^{p,q} (X)\cap W'^{p,q} _2$ be $d$-closed and orthogonal to $ \mathcal{H}^{p,q} _{(2)}(X)$. Then $\alpha$ is $\partial_{J'}\bar{\partial}_{J'}$ -exact:
\[
\alpha=\partial_{J'}\bar{\partial}_{J'}\gamma
\]

for some $\gamma\in W'^{p-1, q-1} _1  (X)\cap C^{p-1, q-1} (X\setminus D)$
\end{lemma}

\begin{proof}
According to the  assumption of being orthogonal to  $ \mathcal{H}^{p,q} _{(2)}(X)$, we must have $\alpha\in \overline{Im \bar{\partial}_{J'}}\oplus  \overline{Im \bar{\partial}^* _{J'}}$.
Since $d= \bar{\partial}_{J'}+ \partial_{J'}$,  and  $\alpha$ is a $d$-closed $(p,q)$-form, it has to be   $ \bar{\partial}_{J'}$-closed, as well. If we set $\alpha=\alpha'+\alpha''$,
where $\alpha'\in \overline{Im \bar{\partial}_{J'}} $ and $\alpha'' \in  \overline{Im \bar{\partial}^* _{J'}}$ , then
\[
\langle \alpha'', \alpha'' \rangle_{g'}=\langle \alpha, \alpha'' \rangle_{g'}= \lim_{n\rightarrow \infty }\langle  \alpha,  \bar{\partial}^* _{J'} \theta_n \rangle_{g'}=\lim_{n\rightarrow \infty }\langle \bar{\partial} _{J'}\alpha, \theta_n \rangle_{g'}=0,
\]

where $\theta_n \in A^{p,q+1} _0 (X\setminus D)$ for $n=1,2,...$ is a sequence such that $\lim_{n\rightarrow \infty} \bar{\partial}^* _{J'} \theta_n =\alpha''$ in $L'^{p,q+1} _2$.

Consequently $\alpha'' =0$ and $\alpha=\alpha'\in  \overline{Im \bar{\partial}_{J'}}$. According to theorem \ref{coh2} there exists 
$\eta\in C^{p,q-1} (X)\cap W'^{p,q-1} _1$ such that $\alpha=\bar{\partial}_{J'} \eta$. Let 
\[
\eta=\eta_0+\eta_1+\eta_2,
\]

where $\eta_0\in \mathcal{H}^{p,q-1}_{(2)}$, $\eta_1\in \overline{Im \partial_{J'}}$ and,   $\eta_2\in \overline{Im \partial^*_{J'}}$. We claim that $\bar{\partial}_{J'}\eta_2 =0$. To prove it we note that there exists a sequence $\{\psi _n\}$ of the elements of   $A^{p+1, q-1} _0 (X\setminus D)$ such that $\partial^* _{J'} \psi_n \rightarrow \eta_2$ in $L'^{p+1,q-1} _2$. Since $\partial _{J'} \alpha =0$ 
and $\partial_{J'}^* \bar{\partial}_{J'}=-\bar{\partial}_{J'}\partial_{J'}^*$ we have
\begin{equation}\label{388}
\begin{split}
0=&\langle \alpha , \partial_{J'}^* \bar{\partial}_{J'}\psi_n  \rangle=\langle   \bar{\partial}_{J'} (\eta_1+\eta_2 ), \partial_{J'}^* \bar{\partial}_{J'}\psi_n   \rangle= \langle   \bar{\partial}_{J'} \eta_1, \partial_{J'}^* \bar{\partial}_{J'}\psi_n   \rangle+\langle   \bar{\partial}_{J'} \eta_2 , \partial_{J'}^* \bar{\partial}_{J'}\psi_n   \rangle.
\end{split}
\end{equation}

We note that 
\begin{equation}\label{399}
 \langle   \bar{\partial}_{J'} \eta_2, \partial_{J'}^* \bar{\partial}_{J'}\psi_n   \rangle\rightarrow  -\langle   \bar{\partial}_{J'} \eta_2,    \bar{\partial}_{J'} \eta_2   \rangle
\end{equation}

and

\begin{equation}\label{400}
 \langle   \bar{\partial}_{J'} \eta_1, \partial_{J'}^* \bar{\partial}_{J'}\psi_n   \rangle =- \langle \partial_{J'} \eta_1,  \bar{\partial}^*_{J'}  \bar{\partial}_{J'}\psi_n   \rangle = 0.
\end{equation}

 From (\ref{388}), (\ref{399}) and (\ref{400} ) it follows that $\bar{\partial}_{J'} \eta_2 =0$. Therefore we can conclude that
$\alpha = \bar{\partial}_{J'} \eta_1$. We also know by theorem (\ref{coh2}) that $\eta_1 =\partial_{J'} \gamma$ for some $\gamma$ 
belonging to $W'^{p-1, q-1} _1  (X)\cap C^{p-1, q-1} (X\setminus D)$. The proof is complete.

\end{proof}

 \section{ Kodaira  theory}\label{sec8}

We start by Kodaira vanishing theorem:

   \subsection{Kodaira Vanishing}

Following our definition of degenerate complex manifolds, we can also study   complex   vector bundels over $X$. The complex vector bundle can also have degenerating complex structure, but for our purpose of developing Kodaira theory the following simple case is sufficient.

%
%
%*****************************************************************
%
%
%
%
%
%end subsubsection de{}.................................................
%........................................................................
%.......................................................................
%
%*************************************************************************
%
%
\begin{defi}\label{linb}
Consider a complex line bundle $\mathcal{L}$ over $X$ which is $J'$ holomorphic on $X\setminus D$.  $\mathcal{L}$ is said  to be $J'$-positive if
its first Chern class contains a $(1,1)$-form $\omega'$ s.t. the
metric $g'$ compatible  with $(\omega ',J')$ is a  
K\"ahler metric in the sense of definition (\ref{def2}).
\end{defi}

%From definition \ref{def2} one can see that the form
%$\gamma$ vanishes along $D$. It then follows that $\mathcal
%{L}|_{D}$, the restriction of the fiber $\mathcal{L}$ to
%$D$,  is trivial and its trivialization extends to a
%neighborhood of $D $.

 The definition of the space     $C^{p,q}
(\mathcal{L})=A^{p,q} (X\setminus D)\cap L'^{p,q} _2 (\mathcal{L}) $ consisting of  smooth $(p,q)$-forms on $X\setminus D$ with values in
$\mathcal{L}$ which are $L'^{p,q} _2$, is straightforward. We equippe  $\mathcal{L}$ with a hermitian metric  such that its curvature becomes equal to  $2\pi i \omega '$. Then we can consider $L'^2$ space of $(p,q)$-forms with values in $\mathcal{L}$. We note that the existence of such a metric follows from  $\partial\bar{\partial }$  lemma  (\ref{ddbr}).  We denote by
$\mathcal{H}^{p,q} _2(X,\mathcal {L})$ the space of $L'^2$ harmonic
$(p,q)$ forms on $X$ with values in $\mathcal{L}$ (defined
similar to section 3.1). One can then show the following
theorem:
\begin{theo}\label{kodaira vanishing}
$\mathcal{H}^{p,q} _{(2)} (X,\mathcal {L})=0 \text { for } p+q>n.$
\end{theo}

\begin{proof}

 We apply $\partial\bar{\partial}$-lemma and corollary \ref{ck}  in order to provide a proof.  Let  the operator $L$ be defined  by
\[
\begin{matrix}
L : C^{p,q}(\mathcal {L}) \rightarrow C^{p,q}(\mathcal{L})\\
L (\eta\otimes s)=\omega' \wedge \eta \otimes s,
\end{matrix}
\]
where $\eta \in C^{p,q}(X)$ and $s \in C^0(\mathcal{L})$. Using
 $\partial \bar{\partial}$-lemma (\ref{ddbr})   the line bundle $\mathcal{L}$ can
be equipped  with a hermitian metric with curvature equal to
$\omega'$. Let
$
\Lambda = L ^*
$
  be the adjoint  of $L$ with respect to $\langle ,\rangle_2 '$. If $D=D'
+D''$ ($D''=\bar \partial_{J'}$) denotes the  connection
associated with the  metric on $\mathcal{L}$, by standard calculation in local coordinates around non-singular points we get the following
relations
\[
[\Lambda,\bar \partial_{J'} ]=-\frac{i}{2} D'^*,
\hspace{0.5cm}
[\Lambda,  L]=(2-p-q)Id.
\]
\\

For the curvature operator  $\Theta$ we have
\[
\Theta \eta= \Theta \wedge \eta =(2\pi/i) L (\eta)=D^2 \eta,
\]
 and
\[
\Theta=D^2=\bar \partial_{J'} D' + D'\bar \partial_{J'}.
\]
Now let $\eta \in \mathcal{H}^{p,q}(\mathcal{L}) $. This implies
that $\bar
\partial_{J'} \eta=0$, and we find

\[
\Theta \eta= \bar \partial_{J'} D' \eta.
\]
and
\[
\begin{split}
2i\langle \Lambda \Theta \eta,\eta\rangle_{g'}  &=  2 i \langle\Lambda \bar
\partial_{J'} D' \eta, \eta\rangle_{g'}= 2i\langle \left(\bar \partial_{J'}
                                   \Lambda -
                                   \frac{i}{2}D'^*\right )D'\eta,\eta\rangle_{g'}\\
\quad                             &=\langle D'^*D' \eta,\eta\rangle_{g'}=\langle D'\eta,D'\eta \rangle_{g'} \geq 0,
\end{split}
\]
where we  have used $\langle \bar \partial_{J'} \Lambda D' \eta,
\eta\rangle_{g'}=\langle \Lambda D'\eta,\bar\partial_{J'}^* \eta\rangle_{g'}=0$. The reason for these identities comes from the fact that 
$\bar{\partial}_{J'} \Lambda D' \eta \in L'^{p,p} _2$ and  this  is because $\eta $ is smooth and $\partial_{J'}$-closed. Moreover according to $\partial\bar{\partial}$  lemma (\ref{ddbr}), the zero-th order part of the connection  belongs to $W' _1$. We can thus apply corollary  (\ref{ck}).   Similarly

\[
\begin{split}
2 i\langle \Theta \Lambda \eta, \eta\rangle_{g'} &= 2i\langle D'\bar \partial_{J'}
\Lambda \eta, \eta\rangle _{g'}=2i\langle D'(\Lambda\bar \partial_{J'} + \frac{i}{2}
    D'^*)\eta,\eta\rangle_{g'}\\
\quad 
    &= -\langle D'D'^* \eta,\eta\rangle_{g'}=-\langle D'^* \eta,D'^* \eta\rangle_{g'} \leq 0.
\end{split}
\]
from which we obtain
\[
2i \langle [\Lambda,\Theta]\eta,\eta\rangle_{g'} \geq 0.
\]
On the other hand $\Theta=(2\pi /i)L$. Thus
\[
2i\langle [\Lambda,\Theta]\eta,\eta \rangle_{g'}   =
4\pi \langle [\Lambda,L]\eta,\eta\rangle_{g'}
= 4\pi(2-p-q) \|\eta\| \geq 0
\]
So for $p+q >2$,  implies that $\eta=0$.
\end{proof}
\

\subsection{Kodaira Embedding Theorem.}  We would like to prove   an
embedding theorem for  $X\setminus U_{\epsilon}$ where $\epsilon
>0 $ is an arbitrary  (small) real number and $U_{\epsilon}$ denotes the $\epsilon$ neighborhood of
$D$. This  theorem seems to be  of importance   due to the following reasons.
Consider the case where $n=2$ and $D$ is a canonical divisor in $X$. According to section \ref{sec2} we can equippe $X$ with a degenerate hyperK\"ahler structure. Assume that $J$ is the initial complex structure and $J'$ denotes the degenerate complex structure on $X$
which is perpendicular to $J$ as described before. 
  By  the argument given in   \cite{Bennequin}, $X\setminus U_{\epsilon}$ is a
pseudo-concave manifold with boundary. It is well known that a generic CR-structure on a 3-dimensional manifold is not embeddable in a projective space (\cite{eps}).  A consequence of our theorem is that the $CR$ structure on $\partial U_{\epsilon}$ is fillable. Moreover
$X\setminus U_{\epsilon}$ admits a compactification which is an ordinary complex manifold.  
  
Also zeros of holomorphic sections of $\mathcal{L}$  in this case give rise to Lagrangian surfaces with respect to $J$.

%bahbah

\begin{theo}\label{kodaira}
Let $(X, D)$ be a singular complex surface in the sense of
definition \ref{def2}. Let $\mathcal{L}$ be a
$J'$-positive line bundle on $X$ in the sense of definition (\ref{linb}). Let $U_{\epsilon}$ be a small
neighborhood of $D$ in $X$ as above. Then, for sufficiently
large values of $n$, holomorphic sections of $\mathcal{L} ^n
|_{X\setminus U_{\epsilon}} $ provide a holomorphic embedding $j:
X\setminus U_{\epsilon} \rightarrow \mathbb{C} \mathbb{P} ^N$ from
$X\setminus U_{\epsilon}$ into some projective space $\mathbb{C}
\mathbb{P} ^N$.
\end{theo}

\begin{proof}  We prove that for large values of $n$ the map
    $j: X\setminus U_{\epsilon} \rightarrow  \mathbb{C}\mathbb{P} ^N$ obtained via holomorphic
   sections of $\mathcal{L} ^n |_{X\setminus U_{\epsilon}}$, separate
    points  and   is regular everywhere.
    Let $p$ and $q$ be two points in $X\setminus
D$. Let $\tilde{X}$ be the degenerate complex manifold obtained from blowing up of $X$ at $p$ and $q$, and
let $\pi: \tilde{X}\rightarrow X $ be the natural projection. We also denote by $\tilde {\mathcal{L}}$ the
pullback of $\mathcal{L}$ with respect to  $\pi$. Set $\pi ^{-1} p = E_p$ and $\pi ^{-1} q = E_q$. Let $\mathcal{U}= \{ u _
\alpha \} _{\alpha\in \mathcal{I}}$ be a special covering of $\tilde{X}$ satisfying the following properties:
\\

i) There exists only one open set $u_ {\alpha _0} \in \mathcal{U}$
such that $u _{\alpha _0} \cap D \neq \emptyset $. Thus we

\hspace{0.5 cm} have $D \subset  u_{\alpha _0}$.
\\

ii) For all $\alpha _{i_1},...\alpha _{i_m} \in \mathcal{I}$ with $m \geq 2$, the intersections $u_{\alpha_{i_1}}
\cap ... \cap u_{\alpha _{i_m}}$  all have

\hspace{0.5 cm} trivial colomologies.
\\

Let $U_p$ and $U_q$ be two small neighborhoods of $E_p$ and $E_q$ respectively.  Let $f_p$ and $f_q$ be two
arbitrary holomorphic sections of $\mathcal{L}^n |_{U_p}$ and $\mathcal{L}^n| _{U_q}$. We would like to show that
there exists a holomorphic section $f$ of $\mathcal{L}^n |_{X\setminus D}$ such that $f|_{E_p} = f_p |_{E_p}$
and $f|_{E_q} = f_q |_{E_q}$. To see this let $ f_{\alpha}: u_{\alpha} \rightarrow \mathcal{L}^n $,  for $\alpha \in \mathcal{I}$,  be a set of holomorphic
sections of $\mathcal{L}^n |_{u_\alpha}$, where $\mathcal{I}$ is an index set.  We suppose that if $u_{\alpha} \cap E_p
  \neq \emptyset $ (resp $u_{\alpha} \cap E_q
  \neq \emptyset $  ) then $f_{\alpha} |_{E_p} = f_p $ (resp
  $f_q$). Now we set $f_{\alpha\beta} = f_{\alpha} - f_{\beta}$ which is defined on
  $u_{\alpha} \cap u_{\beta}$. We are looking for  a set of holomorphic
  sections $g_{\alpha} : u_{\alpha} \rightarrow \mathcal{L}^n$ such that $g_{\alpha} - g_{\beta} =
  f_{\alpha\beta}$ and $g_{\alpha} | _{u_\alpha \cap E_p }= g_{\alpha} | _{u_\alpha \cap E_q} =
  0$, for all $\alpha, \beta$. This is equivalent to say that the collection $\{g_{\alpha}\}$ forms a set of sections of the
  line bundle $\mathcal{L}^n \otimes \mathcal{L}_{p}^{-1} \otimes \mathcal{L}_{q} ^{-1}$, where $\mathcal{L}_p$ (resp
  $\mathcal{L}_q$) are line bundles associated with the divisors $E_p$ (resp.
  $E_q$). For each $\alpha, \beta \in \mathcal{I}$ ,  $ \bar\partial_{J'} f_{\alpha \beta}$
  is a  section of $\Omega ^{0,1} ( \mathcal{L}^n \otimes  \mathcal{L} _{p} ^{-1} \otimes
   \mathcal{L} _{q} ^{-1}) |_{u_{\alpha}
  \cap u_{\beta}}$.  Let $h_{\alpha} \in C^{\infty}( \mathcal{L}^n \otimes  \mathcal{L} _{p} ^{-1}
\otimes \mathcal{L} _{q} ^{-1})
  |_{u_{\alpha}}$ be such that $h_{\alpha} - h_{\beta} =
  f_{\alpha\beta}$.  Then $\{\bar\partial_{J'} h_{\alpha} \}_{\alpha\in \mathcal{I}}$,
    glue together
  to form a global $\bar\partial _{J'}$ -closed section $\theta$ of $A ^{0,1} ( \mathcal{L}^n \otimes
   \mathcal{L} _{p} ^{-1}
   \otimes \mathcal{L} _{q} ^{-1}) 
  $. Now by Kodaira vanishing theorem $\theta$ is also
  $\bar\partial _{J'}$-exact, i.e. $\theta= \bar\partial_{J'} h$ for some $h\in  C^{\infty}( \mathcal{L}^n \otimes
  \mathcal{L} _{p} ^{-1}
  \otimes \mathcal{L}_{q}
  ^{-1})$. This means that $g_{\alpha}= h_{\alpha} - h$ form the
  desired set of holomorphic sections of $\mathcal{L}^n \otimes  \mathcal{L} _{p} ^{-1} \otimes \mathcal{L} _{q}
  ^{-1}$.
\end{proof}
%\begin{cor}
%With the hypothesis of previous theorem, for large values of $k$,  $\mathcal{L}^k$ admits a holomorhic section whose zero section %does not meet $\partial U_{\epsilon}$.
%\end{cor}

%\emph{Proof:} According to \cite{eps} $  \partial U_{\epsilon}$ admits a  Stein filling  $X_{0}$  and thus $(X\setminus U_{\epsilon})%\cup X_0$ is a compactification of $X\setminus U_{\epsilon}$. Obviously  $\mathcal{L}$ can be extended to $\bar{X}$ which is trivial %restricted on $X_0$.
 %\\

\section{Solvability of $\bar{\partial}_{J'}$ in a neighborhood of $D$ } \label{solv}

Consider as before a pair $(X,D)$ consisting of a complex manifold $X$ and a  smooth canonical divisor $D\subset X$. Let $S $ denote the holomorphic section of $K_D=[D]$ vanishing along $D$. Let also $\omega$ denote a K\"ahler form on $X$ representing $c_1 (K_D)$.  We set

\[
S=\alpha+i \beta,
\]
 
where $\alpha, \beta\in \Lambda^2 (X)$. Let $(X, J', \alpha)$ be the degenerate K\"ahler manifold  as defined in section \ref{sec2}.
We   prove the following theorem

\begin{theo}\label{solv9}
In an appropriate neighborhood $U$ of $D$ the equation $\bar{\partial}_{J'} u=f$ admits  a solution $u\in L'^{2} _{(p,q) } (U, loc)$ for every $f\in L'^{2} _{(p,q+1) } (U, loc)$
such that $\bar{\partial}_{J'} f=0$. Moreover if $f\in C^{\infty} \left ( U\setminus D \right )$ then  we have
$u\in C^{\infty} (U\setminus D)$.
\end{theo}

\subsection{The Existence of Plurisubharmonic Maps in the Neighborhood of $D$}\label{secc}

\begin{lemma}\label{lem111} If $\epsilon$ is small enough  and $U_{\epsilon }$ is the neigborhood of $D$ defined in definition (\ref{def2}), then  the map $\phi_0: U_{\epsilon}\rightarrow \mathbb{R}$ defined by $\phi_0 (p) =-|S(p)|^4$,
 where $|.|$ is a hermitian metric on $K_X$, is pluri-subharmonic:
\[
\partial_{J'}\bar{\partial}_{J'} \phi_0 > 0
\] 
\end{lemma}

\begin{proof}
 We consider  a holomorphic coordinates system $(z,w)$ on a neighborhood $U_p$ of a point $p\in D$ such that a $D=\{z=0\}$ and in this coordinates system we take an orthonormal   moving frame $(V_1, V_2)$ given by 

\begin{equation}\label{355}
V_1= \frac{a}{z}\left(  \frac{\partial}{\partial z} -\langle  \frac{\partial}{\partial z},\frac{1}{ \|\frac{\partial }{\partial w}\|^2_{g'}} \frac{\partial}{\partial w}\rangle _{g'}  \frac{\partial}{\partial w}\right), \hspace{1cm} V_2 = b\frac{\partial}{\partial w}.
\end{equation}

where $a$ and $b$ are functions which normalize $V_1$ and $V_2$. Using the canonical coordinates given by  lemma 30 in the appendix A4 of \cite{bah} we can assume that  $da$ and $db$ both vanish at $p$.
This moving frame is $g'$ orthonormal and  belongs to $T'_{\mathbb{C}}X$ (the holomorphic tangent bundle of $X$ with respect to $J$) . Hence the moving frame $(V'_1, V'_2)$ defined by 

\[
V_1 ' := \frac{1}{2}\left ( \Re V_2+ \frac{1}{i} \Re V_1\right ), \hspace{0.5cm} V_2 ' :=  \frac{1}{2}\left ( \Im V_2+ \frac{1}{i} \Im V_1\right ),
\]

is still $g'$-orthonormal and belongs to $T'_{J', \mathbb{C}}X$ the holomorphic tangent space of $X$ with respect to $J'$.  The dual frame
$\{V_1^*, V_2 ^*\}$ is given by 

\begin{equation}\label{33se}
V_1^* = \frac{1}{a}zdz \hspace{1cm} V_2^* =\frac{1}{b} \left (dw + \langle  \frac{\partial}{\partial z}, \frac{1}{\|\frac{\partial }{\partial w}\|^2_{g'}} \frac{\partial}{\partial w}\rangle _{g'}  dz \right ).
\end{equation}

If $z=x+iy$
then $\phi_0 (z,w)=-e^h \left( x^4+y^4+2x^2y^2\right )$, where $e^h$ comes from the hermitian metric on $K_X$. We then have
\[
\partial_{J'}\bar{\partial}_{J'}\phi_0 = \sum _{ij} ^2 \left ( V'_i .\overline{V'_j}. \phi_0 \right ) V'^* _i\wedge  \overline{V'_j}^* +R,
\]
 
where $R$ is obtaine by the terms generated from the action $\bar{\partial}_{J'} V'^* _i$ of $\bar{\partial}_{J'}$  on $V'^* _i$ for $i=1,2$.
It is not difficult to veify that 
\[
\partial_{J'}\bar{\partial}_{J'}\phi_0=\partial_{J'_0}\bar{\partial}_{J'_0}\tilde{\phi}_0+O\left(|\left(x,y\right)|\right).
\]

where $\tilde{\phi_0}= x^4+y^4+2x^2y^2$, and $J'_0$ is the flat degenerate complex structure  described in section (\ref{sec3}).
As we have seen the  operator $\bar{\partial}_{J'}$ is  equal to 

\[
\bar{\partial}_{J'}:= v'_1 \omega' _1+ v'_2 \omega'_2.
\]

where $v'_1, v'_2, w'_1$ and $w'_2$  are defined by  relations (\ref{v})  (\ref{wprim}). Hence the proof of lemma is complete according to  the following computation

\[
v_1 '.\overline{v'_2}  \left(x^4+y^4+2x^2y^2\right ) =-\frac{1}{4}\frac{x\partial_x-y\partial_y}{x^2+y^2}\left (8xy\right ) =0
\]

\[
\begin{split}
v_2 '.\overline{v'_1}  \left(x^4+y^4\right ) &=-\frac{1}{4}\frac{y\partial_x+x\partial_y}{x^2+y^2}\left (  x^2-y^2\right )=0 
\end{split}
\]

\[
\begin{split}
v_1 '.\overline{v'_1}  \left(x^4+y^4\right ) &=-\frac{1}{4}\frac{x\partial_x-y\partial_y}{x^2+y^2}\left (4x^2-4y^2\right ) =-2\\
\end{split}
\]

\[
\begin{split}
v_2 '.\overline{v'_2}  \left(x^4+y^4\right ) &=-\frac{1}{4}\frac{y\partial_x+x\partial_y}{x^2+y^2}\left (8xy \right) =-2\\
\end{split}
\]

\end{proof}

\subsection{The Existence of Moving Frames of $1$ -Forms with Bounded $\bar{\partial}_{J'}$ -Differential}\label{appa}

We take take the same  holomorphic coordinates system $(w,z)$ on a neighborhood $U_p$  of some point $p\in D$ as ini previous subsection. 
 We also define smooth maps $\xi$ and $\eta$ by

\[
\xi+i\eta:=\langle  \frac{\partial}{\partial z},\frac{1}{ \|\frac{\partial }{\partial w}\|^2_{g'}} \frac{\partial}{\partial w}\rangle _{g'}. 
\]

 In the cannonical coordinates constructed   by  lemma 30 in the appendix A4 of \cite{bah}we can assume that $\xi(p)=\eta(p)=0$ .
We also have

\[
\Re V_1 ^* = \frac{1}{a} \left (   xdx-ydy\right ) , \hspace{0.5cm} \Im V_1 ^* = \frac{1}{a}\left (  xdy+ydx \right )
\]
\[
\Re V_2 ^* =  \frac{1}{b}\left (du +\xi dx-\eta dy\right )\hspace{0.5cm} \Im V_2 ^* =\left ( \frac{1}{b}dv+\xi dy+\eta dx \right )
\]

\[
\Re V_1= \frac{1}{2}\left [\frac{ax}{x^2+y^2}\frac{\partial}{\partial x}-\frac{ay}{x^2+y^2}\frac{\partial}{\partial y}+\left(\frac{-\xi x-\eta y}{x^2+y^2} \right) \frac{\partial}{\partial u}+\left (\frac{\xi x-\eta y}{x^2+y^2}\right )\frac{\partial }{\partial v} \right ]
\]
\[
\Im  V_1 =\frac{1}{2}\left [ -\frac{ay}{x^2+y^2}\frac{\partial}{\partial x}-\frac{ax}{x^2+y^2}\frac{\partial}{\partial y}+\left ( \frac{\xi y-\eta x}{x^2+y^2}\right )  \frac{\partial}{\partial u}+\left ( \frac{\xi z-\eta y}{x^2+y^2}\right ) \frac{\partial}{\partial v}\right ]
\]

\[
\Re V_2 = b\frac{\partial}{\partial u} \hspace{0.5cm} \Im V_2 =-b\frac{\partial}{\partial v}
\]

Let us  define $d\bar{\alpha}_i$ and $\frac{\partial}{\partial \bar{\alpha}_i}$ for $i=1,2$ as follows
\[
\begin{split}
d\bar{\alpha}_1 = \Re V_2 ^* -i \Re V_1 ^*  =  \frac{1}{b}\left (du +\xi dx-\eta dy\right )-\frac{i}{a} \left (   xdx-ydy\right )
\end{split}
\]
 
\[
 d\bar{\alpha}_2 =\Im V_2^* +i\Im V_ 1 ^*=\left ( \frac{1}{b}dv+\xi dy+\eta dx \right )+\frac{i}{a}\left (  xdy+ydx \right )
\]
\[
\frac{\partial}{\partial \bar{\alpha}_1}=\frac{1}{2}\left ( \Re V_2 +i \Re V_1\right )=\frac{1}{2}\left (  b\frac{\partial}{\partial u}+ \frac{i}{2}\left [\frac{ax}{x^2+y^2}\frac{\partial}{\partial x}-\frac{ay}{x^2+y^2}\frac{\partial}{\partial y}+\left(\frac{-\xi x-\eta y}{x^2+y^2} \right) \frac{\partial}{\partial u}+\left (\frac{\eta x-\xi y}{x^2+y^2}\right )\frac{\partial }{\partial v} \right ]\right )
\]

\[
\frac{\partial}{\partial \bar{\alpha}_2}=\frac{1}{2}\left ( \Im V_2 -i \Im V_1\right )=\frac{1}{2}\left (  -b\frac{\partial}{\partial v} - \frac{i}{2}\left [ -\frac{ay}{x^2+y^2}\frac{\partial}{\partial x}-\frac{ax}{x^2+y^2}\frac{\partial}{\partial y}+\left ( \frac{\xi y-\eta x}{x^2+y^2}\right )  \frac{\partial}{\partial u}+\left ( \frac{\xi z-\eta y}{x^2+y^2}\right ) \frac{\partial}{\partial v}\right ]\right )
\]

Then we have

\[
\bar{\partial }_{J'}=\frac{\partial}{\partial \bar{\alpha}_1}d\bar{\alpha}_1+ \frac{\partial}{\partial \bar{\alpha}_2}d\bar{\alpha}_2
\]

and we can prove the following lemma

\begin{lemma}\label{lem990}
There exists a constant $C$ such that $\|\bar{\partial}_{J'} \left ( d\bar{\alpha}_i \right )  (p) \|_{g'}\leq C$,  and $C$ is independent of the point $p$.
\end{lemma}

\begin{proof}
The proof follows from a  simple computation. For instance we need to ensure that the terms in 
$\frac{\partial}{\partial \bar{\alpha}_1}\left(     -\frac{i}{a} \left ( xdx-ydy\right )\right )\wedge d\bar{\alpha}_1 $ which  create  singularities  are eliminated by the  terms coming  from the derivatives    
$\frac{\partial}{\partial \bar{\alpha}_2}\left(     -\frac{i}{a} \left ( xdx-ydy\right )\right ) \wedge d\bar{\alpha}_2$. Similarly the singularities created by $\frac{\partial}{\partial \bar{\alpha}_1}\left(   \frac{1}{b}\left (du +\xi dx-\eta dy\right ) \right ) \wedge d\bar{\alpha}_1$
are canceled out by $\frac{\partial}{\partial \bar{\alpha}_2}\left( \frac{1}{b}\left (du +\xi dx-\eta dy\right )  \right ) \wedge d\bar{\alpha}_2$ and vice versa. We are  also using the assumption that $\xi (p)=\eta (p) =0$
\end{proof}

%\textbf{\Large{Appendix}}

%*****************************************************************

\subsection{Proof of Theorem \ref{solv9}}

Let $U$ be a neighborhood of $D$  on which there exists a pluisubharmonic map $\phi_0: U\rightarrow \mathbb{R}$ as considered in lemma (\ref{lem111}). On  the other hand by \cite{Bennequin} (page 9), we know that $ -\log \delta (.,\partial U) $
 satisfies
\[
\partial_{J'}\bar{\partial}_{J'} (-\log \delta (.,\partial U) )>0
\]

in a neighborhood of $\partial U$.  Thus if $\phi_1$ is a smooth function satisfying
\[
\phi_1=\begin{cases}
-\log \delta (.,\partial U)  & \text{in a neighborhood of $\partial U$} \\ 
0 &\text{in a neighborhood of  $D$ }
\end{cases}
\]

 Then  $\lambda\phi_0+ \phi_1$ for large enough $\lambda$  is plurisubharmonic function in $U\setminus D$ and we have
\[
U_c =\{x; x\in U, \phi(x) <c \}\subset\subset U
 \]

We now apply Hormander $L^2$ theory to derive solvability of $\bar{\partial}_{J'}$ in $U$.   According to lemma (\ref{lem990}) in the previous subsecction  near  a given point $p\in U  $, one can choose a $g'$-orthonormal basis of $(1,0)$-forms $\alpha_1, \alpha_2$ in such a way that    $\|\bar\partial _{J'} \alpha_i \|_{g'} <C$  for some constant $C$ independent of $p$.
 Moreover we have
\[
\int  \frac{\partial f}{\partial \alpha^i }  h d\mu_{g'} = \int f\frac{\partial h }{\partial \alpha^i} d\mu_{g'} \hspace{1cm} f,h\in C^{\infty} _0 (U)
\]

Thus by the same  arguments as in \cite{bah2} and (\cite{hor})  one can prove that $D_{p,q+1} (U)$ is dense in $D_{T^*}\cap D_S$ for graph norm 
$f\rightarrow \|f\|_{\phi}+ \|T^* f\|_{\phi}+ \|Sf\|_{\phi}$. Here
\[
T:L'^{2}_{p,q} (U,\phi)\rightarrow L'^2 _{p,q+1} (U,\phi)
\hspace{0.5cm}
S:L'^{2}_{p,q+1} (U,\phi)\rightarrow L'^2 _{p,q+2} (U,\phi)
\]

are  linear closed densely defined operators which extend  $\bar{\partial}_{J'}$.  One then can prove the inequality 
\begin{equation}\label{3990}
\|f\|^2 _{\chi (\phi)} \leq \|T^* f\|^2_{\chi(\phi)} +\|Sf\|_{\chi(\phi)} ^2  \hspace{1cm} f\in D_{(p,q+1)} (U)
\end{equation}
for appropriate constants from  which theorem \ref{solv9} can be deduced.

\appendix

%\section{Proof of theorem 6}\label{appendc}
%\emph{Proof of theorem \ref{kodaira vanishing}:}

\section{Poincar\'e Lemma and DeRham Theory}\label{sec7}

Let  $\mathbb{C}^n =\mathbb{R}^k \times \mathbb{R}^{2n-k}$ be a decomposition of $\mathbb{C}^n$, and let $\pi :\mathbb{C}^n \rightarrow \mathbb{R}^k$ be projection onto the  first component.   Let $g'$ be a bounded   metric on $\mathbb{C}^n \setminus (\mathbb{R}^k\times \{0\} )$. Assume that $g'$ in this local coordinates behaves like part (2) of definition (\ref{def2}). Let $U$ be a contractible neighborhood of the origin, and assume that properties (3) and (4) of definition (\ref{def2}) hold in this neighborhood. Set   $C^m(U)= A^m (U\setminus \mathbb{R}^k\times \{0\}\cap W'^m  _1 (U)$, and as before consider a de-Rham complex

\begin{equation}
...\rightarrow C^{m-1} \xrightarrow{d} C^{m}
\xrightarrow{d} C^{m+1}\rightarrow...
\end{equation}

We denote the cohomology associated with  this complex by $H^m _{(2)} (U)$. Our aim in this section is to show that with the above assumptions  
\[
H^m _{(2)} (U) =0.
\]

Assume that $\alpha\in C^m( U\setminus D)$ satisfies $d\alpha =0$. First, we would like to show that 
  $[\alpha]=0$ in $H^m _{dR} (U\setminus (\{0\}\times \mathbb{C}))$ .
This is obvious for $m\neq 2n-k -1$.  
In dimension $m=2n-k-1$,  we need to show that $\int _{S^{2n-k-1} (\epsilon)}j^* \alpha =0$, where  $S^{2n-k-1} (\epsilon)$ is the $2n-k-1$-sphere of radius $\epsilon$ in $\{0\} \times \mathbb{R}^{2n-k}$, and $\epsilon$ is so small that this sphere is contained in $U$.  Also $j$ is the inclusion map.
But if $\int _{S^{2n-k-1} (\epsilon)} j^*\alpha \neq 0$, we can derive a contradiction by repeating the argument of  lemma (\ref{k}) for the closed $2n-1$ form   $\alpha\wedge \beta$, where $\beta =dx_1\wedge ...dx_k$. 
\\

Now we need to show that the equation $\alpha =d\beta$ has a solution in  $ C^{m-1} (U)$. 
To see this we can proceed as follows. First  write $U$ as a union of some contractible starshaped open sets $\{U_i\}_{i}$, and  apply the standard Poincar\'e lemme in each open set $U_i$.  
\\
 We recall  that if  in  local coordinates $(x_1,...,x_{2n})$ on some $U_i$ a
differential $p$-form $\alpha $ is written as  $\alpha =\sum_I v_{I}
(x) dx_I$, where for $I=(i_1,...,i_p)$, $1\leq i_1 ,..., i_p \leq
2n$, we have $dx_I = dx_{i_1} \wedge ... \wedge dx_{i_p}$. Then  for  $\tilde x \in V_i$

\[
\beta(x)=\sum _{|I|=p, 1\leq k\leq p } (\int_o ^1 t^{p-1} v_I
(t(x-\tilde{x} ) +\tilde{x} ) dt) (-1)^{k-1} (x_{i_k}- \tilde x
_{i_k}) dx_{i_1} \wedge...\widehat{dx_{i_k}}... \wedge dx_{i_p},
\]
 is a solution to $d\beta =\alpha$. It is not difficult to see that if $\alpha \in C^m$ then $\beta\in C^{m-1}$. The only thing that remains is to use a Cech cohomology procedure to correct $\beta$'s  on the intersections of the open sets $U_i$'s to derive a globally defined solution of the equation $d\beta = \alpha$. More precisely, for a good cover $\mathcal{U}=\{U_i\}_{i\in \mathcal{I}}$, we first solve 
$\alpha |_{U_i} =d\beta_i$ for $\beta_i\in C^{m-1} (U_i)  $, $i\in \mathcal{I}$. Then the collection $\{(\beta_i-\beta_j)|_{U_i\cap U_j}\}_{i,j}$ is $d$-closed, and by the above  Poincar\`e homotopy formula we can resolve $(\beta_i-\beta_j)|_{U_i\cap U_j}=d\gamma_{ij}$ ,
for $\gamma_{ij}\in C^{m-2} (U_i\cap U_j )$ $i,j\in\mathcal{I}$. By continuing this procedure,  we  obtain for $p=1,..,m$,  a collection  $\{\beta_{i_1,...,i_p}\}_{i_1,...i_p\in \mathcal{I}}$ of closed $(m-p)$-forms on $\{U_{i_1} \cap...\cap U_{i_{p}}\}_{i_1,...i_p\in \mathcal{I}}$ such that
\[
\sum (-1)^k\beta_{i_1,...\hat{i_k},...i_p}=d\beta_{i_1,...,i_p}.
\]

Here $\beta_{i_1,...,i_{m+1}}$, for all $i_1,...,i_{m+1}$ is  constant on  $U_{i_1}\cap...\cap U_{i_{m+1}}$, and they form an element   $\{\beta_{i_1,...,i_{m+1}}\} \in \delta C^{m}(\mathcal{U}, \mathbb{R})$  ($\delta$-exact in Cech cohomology of constant sheaf). Since we have assumed $[\alpha]=0$ in $H_{dR} ^m$,
we can thus start to correct $\{\beta_{i_{1},...,i_{p}}\}$'s inductively to obtain the correct $\beta$ satisfying $\alpha=d\beta$.
Through this procedure it is clear that $\beta\in C^{m-1} (U)$. Thus we get 

 %bia

\begin{lemma} (Degenerate Poincar\'e Lemma)
Let $X$ be as in the definition (\ref{def2}). Let $U$ be a contractible neighborhood of a point $p\in D$.
Then $H^m _{(2)} (U)=0$.
\end{lemma} 

\begin{cor}\label{isodr}

$H^m _{(2)} (X) \simeq H^m _{dR} (X).$
\end{cor}

Consider the following commutative diagram

\begin{equation}\label{diag}
\begin{tikzcd}[column sep=small, row sep=small] 
    & H^m_{dR} (X) \arrow{ddl}[swap]{j_m} \arrow{ddr}{\eta_{1}} &\\
    & \circlearrowleft &\\[-2ex]
    H^{m} _{(2)} (X) \arrow{rr}[swap]{\eta_2} & &\big ( H^{n-m} _{dR} (X)\big )^* \\
\end{tikzcd}
\end{equation}

Here $j_m$ is  induced  by inclusion, and  $\eta_1$  and $\eta_2$  are respectively the standard  Poincar\'e duality map and the  singular Poincar\`e duality map given by
\[
\eta_2 ([\alpha]) (\beta)=\int_x \alpha\wedge \beta,  \hspace{1cm} [\alpha]_{(2)}\in H^m _{(2)} (X) ,  \hspace{0.5cm} [\beta]\in H^{n-m} _{dR} (X).
\] 

Note that  $\eta_2$ is well defiend due to lemma (\ref{k}).

To see that $j_m$ is an isomorphism, regarding corollary (\ref{isodr}),  it is sufficient to show that $j_m$ is one-to-one. 
If $\alpha\in \Lambda^{m} (X)$ is a closed smooth $m$-form on $X$, and if $[\alpha]_{(2)}=0$ then there exists $\gamma\in C^{m-1}(X) $ such that $\alpha=d\gamma$. It then   follows that $\eta_{2} ( [(\alpha)]_{(2)}  ) =0$ and by diagram (\ref{diag}) $\eta_1 ([\alpha])=0$. This means that 
$[\alpha]=0$ (as an element of $H^m _{dR} (X)$). 

Now take  $[\alpha]\in C^{m}(X)$ s.t $d\alpha =0$. Let $\alpha =\alpha_0 +\alpha_1$ be the $L^2$ decompositon of $\alpha$
into harmonic $\alpha_0 \in \mathcal{H}^m _2$, and $\alpha_1\in \overline{d (A^{m-1} (X\setminus D) \cap L'^m _2)}$. Then since 
$\eta_2 ([\alpha_1]_{(2)}) =0$, it follows that $[\alpha_1]_{(2)} =0$ and thus there exists $\beta_1\in C^{m-1} (X)$ s.t $\alpha_1= d\beta_1$. 
\\

we have also the following theorem as in section 3.1:
\begin{theo}
For each $d$-closed $m$-form $\alpha \in C^{m} (X) $ there exists
a unique  $\alpha_0 \in \mathcal{H} ^{m} _2$ and $\beta \in A^{m-1}
(X)$ s.t. $\alpha =\alpha_0 -d \beta$.  In particular $H^m _{(2)} (X)\simeq \mathcal{H} ^m _2$ 
\end{theo}

 \vspace{1cm}

%\textbf{Statements and Declarations}
%\\

%The authors did not receive support from any organization for the submitted work.\\
%\\

%On behalf of all authors, the corresponding author states that there is no conflict of interest.
%\\

%The authors have no relevant financial or non-financial interests to disclose.

\end{document}